\documentclass[a4paper,9pt]{article}
\usepackage[center]{subfigure}
\usepackage{float}
\usepackage[center,small,bf,ruled]{caption}
\usepackage{vmargin}
\usepackage{amssymb}
\usepackage{amsmath}
\usepackage{graphicx}
\usepackage[latin1]{inputenc}
\usepackage{amsthm}
\usepackage{color}
\usepackage{mathrsfs}
\usepackage{fancyhdr}
\usepackage{setspace}
\usepackage{amsthm}
\usepackage{amssymb}
\usepackage{pstricks, pst-plot, pst-node, pst-tree}
\usepackage{epsfig}
\usepackage{hyperref}
\usepackage{pdfsync}
\usepackage{amsfonts}
\usepackage{graphicx}
\usepackage{makeidx}
\usepackage{txfonts}

\setmarginsrb{35mm}{25mm}{35mm}{25mm}
             {0mm}{10mm}{0mm}{10mm}
\numberwithin{equation}{section}
\makeindex

\begin{document}

\title{Random-time processes governed by differential equations of
fractional distributed order}
\author{L.Beghin\thanks{%
Dep. Statistics, Probability and Appl. Statistics, Sapienza University of
Rome. E-mail address: luisa.beghin@uniroma1.it}}
\date{}
\maketitle

\begin{abstract}
We analyze here different types of fractional differential equations, under
the assumption that their fractional order $\nu \in \left( 0,1\right] $ is
random\ with probability density $n(\nu ).$ We start by considering the
fractional extension of the recursive equation governing the homogeneous
Poisson process $N(t),t>0.$\ We prove that, for a particular (discrete)
choice of $n(\nu )$, it leads to a process with random time, defined as $N(%
\widetilde{\mathcal{T}}_{\nu _{1,}\nu _{2}}(t)),t>0.$ The distribution of
the random time argument $\widetilde{\mathcal{T}}_{\nu _{1,}\nu _{2}}(t)$
can be expressed, for any fixed $t$, in terms of convolutions of
stable-laws. The new process $N(\widetilde{\mathcal{T}}_{\nu _{1,}\nu _{2}})$
is itself a renewal and can be shown to be a Cox process. Moreover we prove
that the survival probability of $N(\widetilde{\mathcal{T}}_{\nu _{1,}\nu
_{2}})$, as well as its probability generating function, are solution to the
so-called fractional relaxation equation of distributed order (see \cite{Vib}%
).

In view of the previous results it is natural to consider diffusion-type
fractional equations of distributed order. We present here an approach to
their solutions in terms of composition of the Brownian motion $B(t),t>0$
with the random time $\widetilde{\mathcal{T}}_{\nu _{1,}\nu _{2}}$. We thus
provide an alternative to the constructions presented in Mainardi and
Pagnini \cite{mapagn} and in Chechkin et al. \cite{che1}, at least in the
double-order case.


\end{abstract}

\section{ Introduction}

In the last decade an increasing attention has been drawn to fractional
extensions of the Poisson process: see, among the others, \cite{re}, \cite%
{Ju}, \cite{Las}, \cite{wa1}, \cite{ma1}, \cite{ma2}, \cite{ma3}, \cite{cao}%
. In particular, the analysis carried out by Beghin and Orsingher \cite{BO}
starts from the generalization of the equation governing the Poisson
process, where the time-derivative is substituted by the fractional
derivative (in the Caputo sense) of order $\nu \in \left( 0,1\right] $:%
\begin{equation}
\frac{d^{\nu }p_{k}}{dt^{\nu }}=-\lambda (p_{k}-p_{k-1}),\quad k\geq 0,
\label{ai3}
\end{equation}%
with initial conditions
\begin{equation}
p_{k}(0)=\left\{
\begin{array}{c}
1\qquad k=0 \\
0\qquad k\geq 1%
\end{array}%
\right.  \label{ai2}
\end{equation}%
and $p_{-1}(t)=0$. The solution to this equation has been expressed as the
density of the random-time process called Fractional Poisson process (FPP)
and defined as%
\begin{equation}
\mathcal{N}_{\nu }(t)=N(\mathcal{T}_{2\nu }(t)),\qquad t>0.  \label{fpp}
\end{equation}%
Here $N$ denotes the standard homogeneous Poisson process with rate
parameter $\lambda >0,$ while $\mathcal{T}_{2\nu }(t),t>0$ is a random
process (independent from $N$) with density given by the folded solution to
the fractional diffusion equation%
\begin{equation}
\frac{\partial ^{2\nu }v}{\partial t^{2\nu }}=c^{2}\frac{\partial ^{2}v}{%
\partial y^{2}},\qquad t>0,\text{ }y\in \mathbb{R}.  \label{fold}
\end{equation}%
Later, in \cite{ba}, the distribution of the FPP $\mathcal{N}_{\nu }$ has
been expressed as

\begin{equation}
p_{k}^{\nu }(t)=(\lambda t^{\nu })^{k}E_{\nu ,\nu k+1}^{k+1}(-\lambda t^{\nu
}),\quad \text{ }k\geq 0,t>0,  \label{gml3}
\end{equation}%
in terms of the so-called Generalized Mittag-Leffler (GML) function, which
is defined as%
\begin{equation}
E_{\alpha ,\beta }^{\gamma }(z)=\sum_{j=0}^{\infty }\frac{\left( \gamma
\right) _{j}\,z^{j}}{j!\Gamma (\alpha j+\beta )},\quad \alpha ,\beta ,\gamma
\in \mathbb{C},\text{ }Re(\alpha ),Re(\beta ),Re(\gamma )>0,  \label{gml2}
\end{equation}%
where $\left( \gamma \right) _{j}=\gamma (\gamma +1)...(\gamma +j-1)$ (for $%
j=1,2,...,$ and $\gamma \neq 0$) is the Pochammer symbol and $\left( \gamma
\right) _{0}=1$ (see \cite{kil}, p.45). Moreover a higher-order
generalization of the previous results has been obtained in \cite{ba} by
introducing \textquotedblleft higher-order fractional
derivatives\textquotedblright\ in (\ref{ai3}) and analyzing the following
equation%
\begin{equation}
\frac{d^{n\nu }p_{k}}{dt^{n\nu }}+\binom{n}{1}\lambda \frac{d^{(n-1)\nu
}p_{k}}{dt^{(n-1)\nu }}+...+\binom{n}{n-1}\lambda ^{n-1}\frac{d^{\nu }p_{k}}{%
dt^{\nu }}=-\lambda ^{n}(p_{k}-p_{k-1}),\quad k\geq 0,  \label{n}
\end{equation}%
where $\nu \in \left( 0,1\right) ,$ subject to the initial conditions%
\begin{eqnarray}
p_{k}(0) &=&\left\{
\begin{array}{c}
1\qquad k=0 \\
0\qquad k\geq 1%
\end{array}%
\right. ,\text{ \quad for }0<\nu <1  \label{n1} \\
\left. \frac{d^{j}}{dt^{j}}p_{k}(t)\right\vert _{t=0} &=&0\qquad
j=1,...,n-1,\quad k\geq 0,\text{ \quad for }\frac{1}{n}<\nu <1  \notag
\end{eqnarray}%
and $p_{-1}(t)=0$. The solution to (\ref{n}) was given by the following
finite sum of GML functions:%
\begin{equation}
\widetilde{p}_{k}^{\nu n}(t)=\sum_{j=1}^{n}\binom{n}{j}\left( \lambda t^{\nu
}\right) ^{n(k+1)-j}E_{\nu ,\nu n(k+1)-\nu j+1}^{kn+n}(-\lambda t^{\nu }).
\label{p}
\end{equation}%
The corresponding process was proved to be a renewal, linked to $\mathcal{N}%
_{\nu }(t),t>0$, by the following relationship%
\begin{equation*}
\widetilde{p}_{k}^{\nu n}(t)=p_{nk}^{\nu }(t)+p_{nk+1}^{\nu
}(t)+...+p_{nk+n-1}^{\nu }(t),\qquad t>0.
\end{equation*}%
Thus it can be interpreted as a FPP which \textquotedblleft
records\textquotedblright\ only the $k$-th order events and disregards the
other ones (for an application to the theory of random motions at finite
velocity, see \cite{BO2}).

We will introduce here the assumption that the fractional order $\nu $ of
the derivative appearing in equation (\ref{ai3}) is itself random, with
distribution $n(\nu ),$ $\nu \in \left( 0,1\right] $: i.e.%
\begin{equation}
\int_{0}^{1}\frac{d^{\nu }p_{k}}{dt^{\nu }}n(\nu )d\nu =-\lambda
(p_{k}-p_{k-1}),\quad k\geq 0,\text{ }\nu \in \left( 0,1\right] .  \label{n2}
\end{equation}%
More precisely, we will concentrate on the case of a double-order discrete
distribution of $\nu $, i.e.%
\begin{equation}
n(\nu )=n_{1}\delta (\nu -\nu _{1})+n_{2}\delta (\nu -\nu _{2}),\qquad 0<\nu
_{1}<\nu _{2}\leq 1,  \label{sax}
\end{equation}%
for $n_{1},n_{2}\geq 0$ and such that $n_{1}+n_{2}=1.$ The assumption (\ref%
{sax}) has been already considered in the context of fractional relaxation
(see \cite{Vib}), as well as for fractional kinetic equations and, in the
last case, it leads to the so-called diffusion with retardation (see \cite%
{mapagn}). As we will see in the next section this assumption on $\nu $
produces a form of the solution which is much more complicated than (\ref%
{gml3}) and (\ref{p}), since it involves infinite sums of GML functions.
Nevertheless the renewal property is still valid and a subordinating
relationship similar to (\ref{fpp}) holds for the corresponding process,
which can be defined as%
\begin{equation}
N(\widetilde{\mathcal{T}}_{\nu _{1},\nu _{2}}(t)),\quad t>0.  \label{p2}
\end{equation}%
In this case the random time $\widetilde{\mathcal{T}}_{\nu _{1},\nu _{2}}$
is represented by a process whose transition density can be expressed either
as an infinite sum of Wright functions or by convolutions of stable laws.

In section 3 we will investigate the relationship between the previous
results and the diffusion equation of fractional distributed order%
\begin{equation}
\int_{0}^{1}\frac{\partial ^{\nu }v}{\partial t^{\nu }}n(\nu )d\nu =\frac{%
\partial ^{2}v}{\partial x^{2}},\quad x\in \mathbb{R},t>0,\;v(x,0)=\delta
(x),  \label{p4}
\end{equation}%
for $0<\nu \leq 1.$\ Equations like (\ref{p4}) have been already studied in
\cite{che1} and \cite{mapagn} in connection with the kinetic description of
anomalous diffusions. It has been proved by Chechkin et al. \cite{che1} that
the solution $v(x,t),x\in \mathbb{R},t>0,$ is a probability density function
and that the corresponding process is subordinated to the Brownian motion
via the following relationship%
\begin{equation}
v(x,t)=\int_{0}^{+\infty }\frac{e^{-x^{2}/4\pi u}}{\sqrt{4\pi u}}G(u,t)du.
\label{sub}
\end{equation}%
In (\ref{sub}) the function $G$ is defined by its Laplace transform%
\begin{equation}
\mathcal{L}\left\{ G(u,t);\eta \right\} =\frac{\int_{0}^{1}n(\nu )\eta ^{\nu
}d\nu }{\eta }e^{-u\int_{0}^{1}n(\nu )\eta ^{\nu }d\nu }.  \label{sub2}
\end{equation}%
In the special case of double-order fractional derivative in (\ref{p4})
these authors focus on the behavior of the second moment of $v(x,t),$ which
suggests that the process can be interpreted as a \textquotedblleft
diffusion with retardation\textquotedblright , in this case. Moreover, under
assumption (\ref{sax}), equation (\ref{p4}) can be seen as a particular case
(for $\gamma =2)$ of the equation (\ref{gu3}) below, which is analyzed in
\cite{sa2}$.$ In this paper only the Fourier transform of the solution is
given in explicit form, in terms of infinite sums of generalized
Mittag-Leffler (GML) functions. Finally the solution to (\ref{p4}) has been
analytically expressed in terms of generalized Wright functions by \cite%
{mapagn}.

We prove here that the solution to equation (\ref{p4}) with the assumption (%
\ref{sax}), i.e.%
\begin{equation}
n_{1}\frac{\partial ^{\nu _{1}}v}{\partial t^{\nu _{1}}}+n_{2}\frac{\partial
^{\nu _{2}}v}{\partial t^{\nu _{2}}}=\frac{\partial ^{2}v}{\partial x^{2}}%
,\quad x\in \mathbb{R},t>0,\;v(x,0)=\delta (x),\;n_{1},n_{2}>0,  \label{aa}
\end{equation}%
for $0<\nu _{1}<\nu _{2}\leq 1$, coincides with the density of the
random-time process%
\begin{equation}
B(\widetilde{\mathcal{T}}_{\nu _{1},\nu _{2}}(t)),\quad t>0,  \label{p3}
\end{equation}%
where $B$ is the standard Brownian motion and the time argument $\widetilde{%
\mathcal{T}}_{\nu _{1},\nu _{2}}$ is the same as in (\ref{p2}), thus writing
in explicit form the density $G$ in (\ref{sub}) as an infinite sum of Wright
functions or by convolutions of stable laws.

Finally we note that the density of the random time $\widetilde{\mathcal{T}}%
_{\nu _{1},\nu _{2}}$ appearing in the processes (\ref{p2}) and (\ref{p3})
coincides with the solution to the equation (\ref{p4}), when a different
hypothesis on the density $n(\nu )$ is assumed, i.e.
\begin{equation}
n(\nu )=n_{1}^{2}\delta (\nu -2\nu _{1})+n_{2}^{2}\delta (\nu -2\nu
_{2})+2n_{1}n_{2}\delta (\nu -(\nu _{1}+\nu _{2})),\qquad 0<\nu _{1}<\nu
_{2}\leq 1,
\end{equation}%
for $n_{1},n_{2}\geq 0$ and such that $n_{1}+n_{2}=1.$ Therefore the
equation governing the process $\widetilde{\mathcal{T}}_{\nu _{1},\nu
_{2}}(t),t>0$ turns out to be%
\begin{equation}
\left( n_{1}\frac{\partial ^{\nu _{1}}v}{\partial t^{\nu _{1}}}+n_{2}\frac{%
\partial ^{\nu _{2}}v}{\partial t^{\nu _{2}}}\right) ^{2}=\frac{\partial
^{2}v}{\partial x^{2}},\quad x\in \mathbb{R},t>0,\;n_{1},n_{2}>0,  \label{bb}
\end{equation}%
for $0<\nu _{1}<\nu _{2}\leq 1$, with the usual initial conditions and $%
v_{t}(x,0)=0$, in addition.

Equations (\ref{aa}) and (\ref{bb}) are proved to govern deeply different
processes: while the former is linked, for any value of $\nu _{1},\nu _{2}$
to a diffusion with retardation (see also \cite{che1} and \cite{che2}), the
same is not true for the second equation, which, depending on the value of
the random indexes, produces a subdiffusion or a superdiffusion.

\section{The recursive equation of distributed order}

\subsection{The double-order fractional case}

We begin by considering the following fractional recursive differential
equation
\begin{equation}
\int_{0}^{1}\frac{d^{\nu }p_{k}}{dt^{\nu }}n(\nu )d\nu =-\lambda
(p_{k}-p_{k-1}),\quad k\geq 0,  \label{uno.2}
\end{equation}%
where, by assumption,
\begin{equation}
n(\nu )\geq 0\text{,\qquad }\int_{0}^{1}n(\nu )d\nu =1,\quad \nu \in \left(
0,1\right]  \label{bi1}
\end{equation}%
subject to the initial conditions%
\begin{equation}
p_{k}(0)=\left\{
\begin{array}{c}
1\qquad k=0 \\
0\qquad k\geq 1%
\end{array}%
\right. ,  \label{uno.1}
\end{equation}%
with $p_{-1}(t)=0$. We apply in (\ref{uno.2}) the definition of fractional
derivative in the sense of Caputo, that is, for $m\in \mathbb{N}$,%
\begin{equation}
\frac{d^{\nu }}{dt^{\nu }}u(t)=\left\{
\begin{array}{l}
\frac{1}{\Gamma (m-\nu )}\int_{0}^{t}\frac{1}{(t-s)^{1+\nu -m}}\frac{d^{m}}{%
ds^{m}}u(s)ds\text{,\qquad for }m-1<\nu <m \\
\frac{d^{m}}{dt^{m}}u(t)\text{,\qquad for }\nu =m%
\end{array}%
\right. ,  \label{dc}
\end{equation}%
(see, for example, \cite{kil}, p.92). As a special case, for $n(\nu )=\delta
(\nu -\overline{\nu }),$ and a particular value of $\overline{\nu }\in
\left( 0,1\right) ,$ equation (\ref{uno.2}) reduces to (\ref{ai3}), which
governs the so-called FPP $\mathcal{N}_{\nu }(t),t>0$ (see, for details,
\cite{ma1}, \cite{ma2} and \cite{ma3}).

In order to get an analytic expression for the solution to (\ref{uno.2}), we
adopt here the following particular form for the density of the fractional
order $\nu $:%
\begin{equation}
n(\nu )=n_{1}\delta (\nu -\nu _{1})+n_{2}\delta (\nu -\nu _{2}),\qquad 0<\nu
_{1}<\nu _{2}\leq 1,  \label{bi3}
\end{equation}%
for $n_{1},n_{2}\geq 0$ and such that $n_{1}+n_{2}=1$ (conditions (\ref{bi1}%
) are trivially fulfilled). The density (\ref{bi3}) has been already used by
\cite{mapagn} and \cite{che1}, in the analysis of the so-called double-order
time-fractional diffusion equation, and corresponds to the case of a
subdiffusion with retardation (see next section for details). Moreover, it
was applied in \cite{Vib} in the context of fractional relaxation with
distributed order.

Under assumption (\ref{bi3}), equation (\ref{uno.2}) becomes%
\begin{equation}
n_{1}\frac{d^{\nu _{1}}p_{k}}{dt^{\nu _{1}}}+n_{2}\frac{d^{\nu _{2}}p_{k}}{%
dt^{\nu _{2}}}=-\lambda (p_{k}-p_{k-1}),\quad k\geq 0.  \label{bi2}
\end{equation}%
By taking the Laplace transform of (\ref{bi2}) we get the following first
result.

\

\noindent \textbf{Theorem 2.1} \ \textit{The Laplace transform of the
solution to equation (\ref{bi2}), under conditions (\ref{uno.1}), is given by%
}%
\begin{equation}
\mathcal{L}\left\{ \widetilde{p}_{k}^{\nu };\eta \right\} =\frac{\lambda
^{k}n_{1}\eta ^{\nu _{1}-1}+\lambda ^{k}n_{2}\eta ^{\nu _{2}-1}}{(\lambda
+n_{1}\eta ^{\nu _{1}}+n_{2}\eta ^{\nu _{2}})^{k+1}},  \label{bi4}
\end{equation}%
\textit{for any }$k\geq 0.$

\noindent \textbf{Proof \ }Formula (\ref{bi4}) can be easily obtained by
applying to (\ref{bi2}) the expression for the Laplace transform of the
Caputo derivative, i.e.%
\begin{eqnarray}
\mathcal{L}\left\{ \frac{d^{\nu }u}{dt^{\nu }};\eta \right\}
&=&\int_{0}^{\infty }e^{-\eta t}\frac{d^{\nu }}{dt^{\nu }}u(t)dt  \label{jen}
\\
&=&\eta ^{\nu }\mathcal{L}\left\{ u(t);\eta \right\} -\sum_{r=0}^{m-1}\eta
^{\nu -r-1}\left. \frac{d^{r}}{dt^{r}}u(t)\right\vert _{t=0},  \notag
\end{eqnarray}%
where $m=\left\lfloor \nu \right\rfloor +1$, which yields, for $k\geq 1,$%
\begin{equation}
n_{1}\eta ^{\nu _{1}}\mathcal{L}\left\{ \widetilde{p}_{k}^{\nu };\eta
\right\} +n_{2}\eta ^{\nu _{2}}\mathcal{L}\left\{ \widetilde{p}_{k}^{\nu
};\eta \right\} =-\lambda \left[ \mathcal{L}\left\{ \widetilde{p}_{k}^{\nu
};\eta \right\} -\mathcal{L}\left\{ \widetilde{p}_{k-1}^{\nu };\eta \right\} %
\right] .  \label{bi5}
\end{equation}%
By recursively using (\ref{bi5}) we get%
\begin{equation}
\mathcal{L}\left\{ \widetilde{p}_{k}^{\nu };\eta \right\} =\left( \frac{%
\lambda }{\lambda +n_{1}\eta ^{\nu _{1}}+n_{2}\eta ^{\nu _{2}}}\right) ^{k}%
\mathcal{L}\left\{ \widetilde{p}_{0}^{\nu };\eta \right\} ,\;k\geq 1.
\label{bi6}
\end{equation}%
For $k=0,$ we get, instead:%
\begin{equation}
\mathcal{L}\left\{ \widetilde{p}_{0}^{\nu };\eta \right\} =\frac{n_{1}\eta
^{\nu _{1}-1}+n_{2}\eta ^{\nu _{2}-1}}{\lambda +n_{1}\eta ^{\nu
_{1}}+n_{2}\eta ^{\nu _{2}}},  \label{lo}
\end{equation}%
which, together with (\ref{bi6}), gives (\ref{bi4}).\hfill $\square $

\

The Laplace transform (\ref{bi4}) can be compared with formula (4.8) of \cite%
{mee}, where a Poisson process time-chenged by an arbitrary subordinator is
considered.

We can not use a direct method in order to invert analytically the Laplace
transform (\ref{bi4}). Indeed an explicit inversion formula is available
only for $k=0,$ while for $k>0$ the presence of the power $k+1$ makes the
analytic inversion too complicated. For $k=0,$ we can apply the well-known
expression of the Laplace transform of the GML function defined in (\ref%
{gml2}) (see \cite{kil}, p.47), i.e.%
\begin{equation}
\mathcal{L}\left\{ t^{\gamma -1}E_{\beta ,\gamma }^{\delta }(\omega t^{\beta
});\eta \right\} =\frac{\eta ^{\beta \delta -\gamma }}{(\eta ^{\beta
}-\omega )^{\delta }},  \label{pra}
\end{equation}%
(where $\text{Re}(\beta )>0,$ $\text{Re}(\gamma )>0,$ $\text{Re}(\delta )>0$
and $s>|\omega |^{\frac{1}{\text{Re}(\beta )}})$ and the resulting formulae
(26) and (27) of \cite{sa2}.

Therefore we get%
\begin{eqnarray}
&&\widetilde{p}_{0}^{\nu }(t)  \label{pon} \\
&=&\sum_{r=0}^{\infty }\left( -\frac{n_{1}t^{\nu _{2}-\nu _{1}}}{n_{2}}%
\right) ^{r}E_{\nu _{2},(\nu _{2}-\nu _{1})r+1}^{r+1}\left( -\frac{\lambda
t^{\nu _{2}}}{n_{2}}\right) -\sum_{r=0}^{\infty }\left( -\frac{n_{1}t^{\nu
_{2}-\nu _{1}}}{n_{2}}\right) ^{r+1}E_{\nu _{2},(\nu _{2}-\nu
_{1})(r+1)+1}^{r+1}\left( -\frac{\lambda t^{\nu _{2}}}{n_{2}}\right) ,
\notag
\end{eqnarray}%
under condition $|n_{1}\eta ^{\nu _{1}}/(n_{2}\eta ^{\nu _{2}}+\lambda )|<1$
(which is fulfilled, for $\nu _{2}>\nu _{1},$ $\lambda >0).$

For $k>0,$ we adopt an approach similar to those used in \cite{be}, \cite{BO}%
, \cite{ob1}, \cite{ob2}, \cite{orpo} and \cite{orpo2} (for different types
of fractional differential equations), which leads to an expression of the
solution in terms of convolutions of known distributions. In particular we
will resort to the class of completely asymmetric stable laws (of index less
than one). More precisely, let us denote by $\overline{p}_{\alpha }(\cdot
;z) $, for $j=1,2$, the density of a stable random variable $\mathcal{X}%
_{\alpha }$ of index $\alpha \in \left( 0,1\right) $ and parameters equal to
$\beta =1,$ $\mu =0$ and $\sigma =\left( |z|\cos \frac{\pi \alpha }{2}%
\right) ^{1/\alpha }$ (see \cite{sam} for the definitions and the properties
of this class of stable laws). As well-known, $\mathcal{X}_{\alpha }$ is
endowed by the following Laplace transform%
\begin{equation}
\mathcal{L}\left\{ \overline{p}_{\alpha }(\cdot ;y);\eta \right\} =e^{-\frac{%
\sigma ^{\alpha }}{\cos (\pi \alpha /2)}\eta ^{\alpha }},  \label{sta}
\end{equation}%
which will be particularly useful in inverting (\ref{bi4}). We need moreover
the following result proved in \cite{ob1}: the solution to the following
fractional diffusion equation%
\begin{equation}
\left\{
\begin{array}{l}
\frac{\partial ^{2\alpha }v}{\partial t^{2\alpha }}=c^{2}\frac{\partial ^{2}v%
}{\partial y^{2}},\qquad t>0,\text{ }y\in \mathbb{R}\text{, }c\in \mathbb{R}
\\
v(y,0)=\delta (y),\qquad \text{for }0<\alpha <1 \\
v_{t}(y,0)=0,\qquad \text{for }1/2<\alpha <1%
\end{array}%
\right. .  \label{due.16}
\end{equation}%
can be expressed as%
\begin{equation}
v_{2\alpha }(y,t)=\frac{1}{2c\Gamma (1-\alpha )}\int_{0}^{t}\frac{\overline{p%
}_{\alpha }(s;|y|)}{(t-s)^{\alpha }}ds=\frac{1}{2c}I^{\alpha }\left\{
\overline{p}_{\alpha }(\cdot ;|y|)\right\} (t),\qquad t>0,\text{ }y\in
\mathbb{R},  \label{due.16b}
\end{equation}%
where $I^{\alpha }\left\{ \cdot \right\} $ denotes the Riemann-Liouville
fractional integral of order $\alpha .$\ By $\overline{v}_{2\alpha }(y,t)$
we will denote the folded solution to (\ref{due.16}), i.e.%
\begin{equation}
\overline{v}_{2\alpha }(y,t)=\left\{
\begin{array}{l}
2v_{2\alpha }(y,t),\qquad \text{ }y>0 \\
0,\qquad y<0%
\end{array}%
\right. .  \label{due.16f}
\end{equation}

\

\noindent \textbf{Theorem 2.2} \ \textit{The solution to equation (\ref{bi2}%
), under conditions (\ref{uno.1}), is given, for any }$k\geq 0,$\textit{\
and }$t>0,$\textit{\ by}%
\begin{eqnarray}
&&\widetilde{p}_{k}^{\nu }(t)=\int_{0}^{+\infty }p_{k}(y)q_{\nu _{1},\nu
_{2}}(y,t)dy  \label{bi7} \\
&=&\frac{1}{k!}\int_{0}^{+\infty }y^{k}e^{-y}\left[ \int_{0}^{t}\overline{p}%
_{\nu _{2}}(t-s;y)\overline{v}_{2\nu _{1}}(y,s)ds+\int_{0}^{t}\overline{p}%
_{\nu _{1}}(t-s;y)\overline{v}_{2\nu _{2}}(y,s)ds\right] dy,  \notag
\end{eqnarray}%
\textit{where }$p_{k},$\textit{\ }$k\geq 0$\textit{, represents the
distribution of the standard homogeneous Poisson process }$N(t),t>0$\textit{%
\ (with intensity }$1$) \textit{and }$\overline{p}_{\nu _{j}}(\cdot ;z)$%
\textit{\ denotes the density of the stable random variable }$X_{\nu _{j}}$%
\textit{\ of index }$\nu _{j}\in \left( 0,1\right) ,$\textit{\ for }$j=1,2,$%
\textit{\ with parameters equal }$\beta =1,$\textit{\ }$\mu =0$\textit{\ and
}$\sigma =\left( \frac{n_{j}}{\lambda }|y|\cos \frac{\pi \nu _{j}}{2}\right)
^{1/\nu _{j}}$\textit{.}

\noindent \textbf{Proof \ }We observe that (\ref{bi4}) can be rewritten as
follows%
\begin{eqnarray}
&&\mathcal{L}\left\{ \widetilde{p}_{k}^{\nu };\eta \right\} =(\frac{n_{1}}{%
\lambda }\eta ^{\nu _{1}-1}+\frac{n_{2}}{\lambda }\eta ^{\nu _{2}-1})%
\mathcal{L}\left\{ p_{k};\frac{n_{1}}{\lambda }\eta ^{\nu _{1}}+\frac{n_{2}}{%
\lambda }\eta ^{\nu _{2}}\right\}  \label{bi8} \\
&=&(\frac{n_{1}}{\lambda }\eta ^{\nu _{1}-1}+\frac{n_{2}}{\lambda }\eta
^{\nu _{2}-1})\int_{0}^{+\infty }p_{k}(y)e^{-(\frac{n_{1}}{\lambda }\eta
^{\nu _{1}}+\frac{n_{2}}{\lambda }\eta ^{\nu _{2}})y}dy,  \notag
\end{eqnarray}%
since, for the distribution of $N$, the Laplace transform reads%
\begin{equation*}
\mathcal{L}\left\{ p_{k};\eta \right\} =\frac{1}{(1+\eta )^{k+1}},\quad
k\geq 0.
\end{equation*}%
The exponential in (\ref{bi8}) coincides with the Laplace transform of the
following convolution of the stable laws $\overline{p}_{\nu _{1}}$ and $%
\overline{p}_{\nu _{2}}$:%
\begin{equation}
g_{\nu _{1},\nu _{2}}(w;y)=\int_{0}^{w}\overline{p}_{\nu _{1}}(w-x;y)%
\overline{p}_{\nu _{2}}(x;y)dx.  \label{bi9}
\end{equation}%
Therefore, by considering that%
\begin{equation*}
\eta ^{\nu -1}=\frac{1}{\Gamma (1-\nu )}\int_{0}^{+\infty }e^{-\eta
t}t^{-\nu }dt,
\end{equation*}%
we obtain%
\begin{eqnarray}
\widetilde{p}_{k}^{\nu }(t) &=&\frac{n_{1}}{\lambda \Gamma (1-\nu _{1})}%
\int_{0}^{t}(t-w)^{-\nu _{1}}\mathcal{L}^{-1}\left\{ \int_{0}^{+\infty
}p_{k}(y)e^{-(\frac{n_{1}}{\lambda }\eta ^{\nu _{1}}+\frac{n_{2}}{\lambda }%
\eta ^{\nu _{2}})y}dy;w\right\} dw+  \notag \\
&&+\frac{n_{2}}{\lambda \Gamma (1-\nu _{2})}\int_{0}^{t}(t-w)^{-\nu _{2}}%
\mathcal{L}^{-1}\left\{ \int_{0}^{+\infty }p_{k}(y)e^{-(\frac{n_{1}}{\lambda
}\eta ^{\nu _{1}}+\frac{n_{2}}{\lambda }\eta ^{\nu _{2}})y}dy;w\right\} dw
\notag \\
&=&\frac{n_{1}}{\lambda \Gamma (1-\nu _{1})}\int_{0}^{t}(t-w)^{-\nu
_{1}}\left( \int_{0}^{+\infty }p_{k}(y)g_{\nu _{1},\nu _{2}}(w;y)dy\right)
dw+  \label{bi} \\
&&+\frac{n_{2}}{\lambda \Gamma (1-\nu _{2})}\int_{0}^{t}(t-w)^{-\nu
_{2}}\left( \int_{0}^{+\infty }p_{k}(y)g_{\nu _{1},\nu _{2}}(w;y)dy\right)
dw.  \notag
\end{eqnarray}%
By inserting (\ref{bi9}) into (\ref{bi}) and changing the integration's
order, we get%
\begin{eqnarray*}
\widetilde{p}_{k}^{\nu }(t) &=&\frac{n_{1}}{\lambda \Gamma (1-\nu _{1})}%
\int_{0}^{+\infty }p_{k}(y)\left( \int_{0}^{t}\overline{p}_{\nu
_{2}}(x;y)dx\int_{x}^{t}(t-w)^{-\nu _{1}}\overline{p}_{\nu
_{1}}(w-x;y)dw\right) dy \\
&&+\frac{n_{2}}{\lambda \Gamma (1-\nu _{2})}\int_{0}^{+\infty
}p_{k}(y)\left( \int_{0}^{t}\overline{p}_{\nu
_{1}}(x;y)dx\int_{x}^{t}(t-w)^{-\nu _{2}}\overline{p}_{\nu
_{2}}(w-x;y)dw\right) dy \\
&=&\frac{n_{1}}{\lambda }\int_{0}^{+\infty }p_{k}(y)\left( \int_{0}^{t}%
\overline{p}_{\nu _{2}}(x;y)I^{\nu _{1}}\left\{ \overline{p}_{\nu
_{1}}(\cdot ;y)\right\} (t-x)dx\right) dy \\
&&+\frac{n_{2}}{\lambda }\int_{0}^{+\infty }p_{k}(y)\left( \int_{0}^{t}%
\overline{p}_{\nu _{1}}(x;y)I^{\nu _{2}}\left\{ \overline{p}_{\nu
_{2}}(\cdot ;y)\right\} (t-x)dx\right) dy.
\end{eqnarray*}%
By considering (\ref{due.16b}) and (\ref{due.16f}), for $c=\lambda /n_{j},$
for $j=1,2,$ formula (\ref{bi7}) immediately follows.\hfill $\square $

\

\noindent \textbf{Remark 2.1 \ }The previous result shows that the solution
to (\ref{bi2}) can be expressed as the probability distribution of a
standard Poisson process $N(t),t>0,$ composed with a random time argument
with transition density $q_{\nu _{1},\nu _{2}}(y,t)$, that will be denoted
as $\widetilde{\mathcal{T}}_{\nu _{1},\nu _{2}}$ (independent from $N$):
thus we can write%
\begin{equation}
\widetilde{p}_{k}^{\nu }(t)=\Pr \left\{ N(\widetilde{\mathcal{T}}_{\nu
_{1},\nu _{2}}(t))=k\right\} ,\qquad k\geq 0,\text{ }t>0.  \label{bi10}
\end{equation}%
It is proved in \cite{BO} that the solution to the fractional equation (\ref%
{ai3}) is the density of the composition of $N(t),t>0$ with a random time
argument $\mathcal{T}_{\nu }(t)$, whose density is given by $\overline{v}%
_{2\nu }(y,t).$ The properties of this process have been extensively
analyzed in \cite{ba}: it turns out to be a Cox process, with directing
measure equal to $\Lambda \left( \left( 0,t\right] \right) \equiv \mathcal{T}%
_{\nu }(t).$ We will prove below that an analogous result is valid for the
process $\widetilde{\mathcal{N}}_{\nu _{1},\nu _{2}}(t)\equiv N(\widetilde{%
\mathcal{T}}_{\nu _{1},\nu _{2}}(t))$ introduced here. Moreover we will
check that it is also a renewal process.

\

We derive now a series expression for the transition density $q_{\nu
_{1},\nu _{2}}(y,t)$ of the random time-argument $\widetilde{\mathcal{T}}%
_{\nu _{1},\nu _{2}}(t),t>0$, which is alternative to the integral one given
in Theorem 2.2.

\

\noindent \textbf{Theorem 2.3} \ \textit{The density }$q_{\nu _{1},\nu
_{2}}(y,t)$\textit{\ of the random time-argument }$\widetilde{\mathcal{T}}%
_{\nu _{1},\nu _{2}}(t),t>0$\textit{\ can be expressed as follows}%
\begin{eqnarray}
&&q_{\nu _{1},\nu _{2}}(y,t)  \label{bi11} \\
&=&\frac{n_{1}}{\lambda t^{\nu _{1}}}\sum_{r=0}^{\infty }\frac{1}{r!}\left( -%
\frac{n_{2}|y|}{\lambda t^{\nu _{2}}}\right) ^{r}\mathcal{W}_{-\nu
_{1},1-\nu _{2}r-\nu _{1}}\left( -\frac{n_{1}|y|}{\lambda t^{\nu _{1}}}%
\right) +\frac{n_{2}}{\lambda t^{\nu _{2}}}\sum_{r=0}^{\infty }\frac{1}{r!}%
\left( -\frac{n_{1}|y|}{\lambda t^{\nu _{1}}}\right) ^{r}\mathcal{W}_{-\nu
_{2},1-\nu _{1}r-\nu _{2}}\left( -\frac{n_{2}|y|}{\lambda t^{\nu _{2}}}%
\right) ,  \notag
\end{eqnarray}%
\textit{where}%
\begin{equation*}
\mathcal{W}_{\alpha ,\beta }(x)=\sum_{k=0}^{\infty }\frac{x^{k}}{k!\Gamma
(\alpha k+\beta )},\qquad \alpha >-1,\text{ }\beta \in \mathbb{C},\text{ }%
x\in \mathbb{R}
\end{equation*}%
\textit{is the Wright function.}

\noindent \textbf{Proof \ }We recall that the solution to the diffusion
equation (\ref{due.16}) can be expressed as%
\begin{equation}
v_{2\alpha }(y,t)=\frac{1}{2ct^{\alpha }}\mathcal{W}_{-\alpha ,1-\alpha
}\left( -\frac{|y|}{ct^{\alpha }}\right) ,\qquad t>0,\text{ }y\in \mathbb{R}%
\text{,}  \label{bi12}
\end{equation}%
(see \cite{mai}, for details). Then we get from (\ref{bi7}) that
\begin{eqnarray}
&&q_{\nu _{1},\nu _{2}}(y,t)  \label{fel3} \\
&=&\frac{n_{1}}{\lambda }\int_{0}^{t}\overline{p}_{\nu _{2}}(t-s;|y|)\frac{1%
}{s^{\nu _{1}}}\mathcal{W}_{-\nu _{1},1-\nu _{1}}\left( -\frac{n_{1}|y|}{%
\lambda s^{\nu _{1}}}\right) ds+\frac{n_{2}}{\lambda }\int_{0}^{t}\overline{p%
}_{\nu _{1}}(t-s;|y|)\frac{1}{s^{\nu _{2}}}\mathcal{W}_{-\nu _{2},1-\nu
_{2}}\left( -\frac{n_{2}|y|}{\lambda s^{\nu _{2}}}\right) ds.  \notag
\end{eqnarray}%
We now consider the series representation of the stable law of order $\alpha
\in \left( 0,1\right) $ given in \cite{fel} (formula (6.10), p.583) and
already used (with some corrections), in the fractional context, in \cite%
{ob2}:
\begin{equation}
\overline{p}_{\alpha }(x;\gamma ,1)=\frac{\alpha }{\pi }\sum_{r=0}^{\infty
}(-1)^{r}\frac{\Gamma (\alpha (r+1))}{r!}x^{-\alpha (r+1)-1}\sin \left[
\frac{\pi }{2}(\gamma +\alpha )(r+1)\right] .  \label{fel}
\end{equation}%
In (\ref{fel}) the canonical Feller representation for the stable laws (with
null position parameter $\mu $) has been used, i.e.%
\begin{equation*}
\overline{p}_{\alpha }(x;\gamma ,\zeta )=\frac{1}{2\pi }\int_{-\infty
}^{+\infty }e^{-i\theta x}\exp \left\{ -\zeta |\theta |^{\alpha }e^{-\frac{%
i\pi }{2}\frac{\theta }{|\theta |}}\right\} d\theta ,\qquad \alpha \neq 1;
\end{equation*}%
hence we must convert the parameters appearing there into those used here,
as follows:%
\begin{eqnarray*}
\alpha &=&\nu _{j} \\
\gamma &=&\frac{2}{\pi }\arctan \left( \tan \frac{\pi \nu _{j}}{2}\right)
=\nu _{j}=\alpha \\
\zeta &=&\frac{\sigma ^{\nu _{j}}}{\cos \frac{\pi \nu _{j}}{2}}=\frac{%
n_{j}|y|\cos \frac{\pi \nu _{j}}{2}}{\lambda \cos \frac{\pi \nu _{j}}{2}}=%
\frac{n_{j}}{\lambda }|y|.
\end{eqnarray*}%
By taking into account the self-similarity property, the stable densities
appearing in (\ref{fel3}) become%
\begin{eqnarray}
\overline{p}_{\nu _{j}}(x;\nu _{j},\frac{n_{j}}{\lambda }|y|) &=&\frac{%
\lambda ^{1/\nu _{j}}}{\left( n_{j}|y|\right) ^{1/\nu _{j}}}\overline{p}%
_{\nu _{j}}\left( \frac{x\lambda ^{1/\nu _{j}}}{\left( n_{j}|y|\right)
^{1/\nu _{j}}};\nu _{j},1\right)  \label{fel2} \\
&=&\frac{\lambda ^{1/\nu _{j}}}{\pi \left( n_{j}|y|\right) ^{1/\nu _{j}}}%
\sum_{r=0}^{\infty }(-1)^{r-1}\frac{\Gamma (\nu _{j}r+1)}{r!}\left( \frac{%
x\lambda ^{1/\nu _{j}}}{\left( n_{j}|y|\right) ^{1/\nu _{j}}}\right) ^{-\nu
_{j}r-1}\sin (\pi \nu _{j}r)  \notag \\
&=&\frac{1}{\pi x}\sum_{r=0}^{\infty }(-1)^{r-1}\frac{\Gamma (\nu _{j}r+1)}{%
r!}\left( \frac{n_{j}|y|}{\lambda ^{1/\nu _{j}}x^{\nu _{j}}}\right) ^{r}\sin
(\pi \nu _{j}r),  \notag
\end{eqnarray}%
so that (\ref{fel3}) reads%
\begin{eqnarray}
&&q_{\nu _{1},\nu _{2}}(y,t)  \label{fel4} \\
&=&\frac{n_{1}}{\lambda \pi }\sum_{r=0}^{\infty }(-1)^{r-1}\frac{\Gamma (\nu
_{2}r+1)}{r!}\left( \frac{n_{2}}{\lambda }|y|\right) ^{r}\sin \left( \pi \nu
_{2}r\right) \int_{0}^{t}\frac{1}{(t-s)^{\nu _{2}r+1}s^{\nu _{1}}}\mathcal{W}%
_{-\nu _{1},1-\nu _{1}}\left( -\frac{n_{1}|y|}{\lambda s^{\nu _{1}}}\right)
ds  \notag \\
&&+\frac{n_{2}}{\lambda \pi }\sum_{r=0}^{\infty }(-1)^{r-1}\frac{\Gamma (\nu
_{1}r+1)}{r!}\left( \frac{n_{1}}{\lambda }|y|\right) ^{r}\sin \left( \pi \nu
_{1}r\right) \int_{0}^{t}\frac{1}{(t-s)^{\nu _{1}r+1}s^{\nu _{2}}}\mathcal{W}%
_{-\nu _{2},1-\nu _{2}}\left( -\frac{n_{2}|y|}{\lambda s^{\nu _{2}}}\right)
ds  \notag \\
&=&\frac{n_{1}}{\lambda \pi }\sum_{r=0}^{\infty }(-1)^{r-1}\frac{\Gamma (\nu
_{2}r+1)}{r!}\left( \frac{n_{2}}{\lambda }|y|\right) ^{r}\sin \left( \pi \nu
_{2}r\right) \sum_{l=0}^{\infty }\frac{\left( -n_{1}|y|/\lambda \right) ^{l}%
}{l!\Gamma \left( -\nu _{1}l+1-\nu _{1}\right) }\int_{0}^{t}\frac{1}{%
(t-s)^{\nu _{2}r+1}s^{\nu _{1}(l+1)}}ds+  \notag \\
&&+\frac{n_{2}}{\lambda \pi }\sum_{r=0}^{\infty }(-1)^{r-1}\frac{\Gamma (\nu
_{1}r+1)}{r!}\left( \frac{n_{1}}{\lambda }|y|\right) ^{r}\sin (\pi \nu
_{1}r)\sum_{l=0}^{\infty }\frac{\left( -n_{2}|y|/\lambda \right) ^{l}}{%
l!\Gamma \left( -\nu _{2}l+1-\nu _{2}\right) }\int_{0}^{t}\frac{1}{%
(t-s)^{\nu _{1}r+1}s^{\nu _{2}(l+1)}}ds  \notag \\
&=&\frac{n_{1}}{\lambda \pi t^{\nu _{1}}}\sum_{r=0}^{\infty }(-1)^{r-1}\frac{%
\Gamma (\nu _{2}r+1)}{r!}\left( \frac{n_{2}|y|}{\lambda t^{\nu _{2}}}\right)
^{r}\sin \left( \pi \nu _{2}r\right) \sum_{l=0}^{\infty }\frac{1}{l!}\frac{%
\Gamma (-\nu _{2}r)}{\Gamma \left( 1-\nu _{1}l-\nu _{1}-\nu _{2}r\right) }%
\left( -\frac{n_{1}|y|}{\lambda t^{\nu _{1}}}\right) ^{l}+  \notag \\
&&+\frac{n_{2}}{\lambda \pi t^{\nu _{2}}}\sum_{r=0}^{\infty }(-1)^{r-1}\frac{%
\Gamma (\nu _{1}r+1)}{r!}\left( \frac{n_{1}|y|}{\lambda t^{\nu _{2}}}\right)
^{r}\sin \left( \pi \nu _{1}r\right) \sum_{l=0}^{\infty }\frac{1}{l!}\frac{%
\Gamma (-\nu _{1}r)}{\Gamma \left( 1-\nu _{2}l-\nu _{2}-\nu _{1}r\right) }%
\left( -\frac{n_{2}|y|}{\lambda t^{\nu _{2}}}\right) ^{l}  \notag \\
&=&\frac{n_{1}}{\lambda t^{\nu _{1}}}\sum_{r=0}^{\infty }\frac{1}{r!}\left( -%
\frac{n_{2}|y|}{\lambda t^{\nu _{2}}}\right) ^{r}\sum_{l=0}^{\infty }\frac{1%
}{l!\Gamma \left( 1-\nu _{1}l-\nu _{1}-\nu _{2}r\right) }\left( -\frac{%
n_{1}|y|}{\lambda t^{\nu _{1}}}\right) ^{l}+  \notag \\
&&+\frac{n_{2}}{\lambda t^{\nu _{2}}}\sum_{r=0}^{\infty }\frac{1}{r!}\left( -%
\frac{n_{1}|y|}{\lambda t^{\nu _{1}}}\right) ^{r}\sum_{l=0}^{\infty }\frac{1%
}{l!\Gamma \left( 1-\nu _{2}l-\nu _{2}-\nu _{1}r\right) }\left( -\frac{%
n_{2}|y|}{\lambda t^{\nu _{2}}}\right) ^{l},  \notag
\end{eqnarray}%
where, in the last step, we have used the reflection formula of the Gamma
function.\hfill $\square $

\

\noindent \textbf{Remark 2.2 \ }Consider the special case $n_{1}=0$, $%
n_{2}=1 $: the distribution of the random order $\nu $ reduces, in this
case, to%
\begin{equation}
n(\nu )=\delta (\nu -\nu _{2}),\qquad 0<\nu _{2}\leq 1,
\end{equation}%
so that equation (\ref{bi2}) becomes the fractional equation (\ref{ai3}),
with $\nu =\nu _{2}$. The Laplace transform (\ref{bi4}) simplifies to
\begin{equation}
\mathcal{L}\left\{ \widetilde{p}_{k}^{\nu };\eta \right\} =\frac{\lambda
^{k}\eta ^{\nu _{2}-1}}{(\lambda +\eta ^{\nu _{2}})^{k+1}},
\end{equation}%
thus giving the well-known distribution%
\begin{equation*}
p_{k}^{\nu }(t)=\lambda ^{k}t^{\nu _{2}}E_{\nu _{2},\nu _{2}k+1}^{k+1}\left(
-\lambda t^{\nu _{2}}\right) ,
\end{equation*}%
for any $k\geq 0$ (see \cite{ba} for details). The result of Theorem 2.2 can
be specialized as follows, for $n_{1}=0$, $n_{2}=1$:%
\begin{equation}
p_{k}^{\nu }(t)=\int_{0}^{+\infty }p_{k}(y)q_{\nu _{1},\nu _{2}}(y,t)dy=%
\frac{1}{k!}\int_{0}^{+\infty }y^{k}e^{-y}\overline{v}_{2\nu _{2}}(y,t)dy.
\label{spe}
\end{equation}%
In (\ref{spe}) we have taken into account that the density of the stable
random variable $\mathcal{X}_{\nu _{1}}$ with $\mu =0$ and $\sigma =\left(
\frac{n_{1}}{\lambda }y\cos \frac{\pi \nu _{1}}{2}\right) ^{1/\nu _{1}}$
degenerates to the Dirac's delta function (i.e. $\overline{p}_{\nu
_{1}}(w-x;y)=\delta (x-w)$)$,$ so that the density (\ref{bi9}) becomes $%
g_{\nu _{2}}(w;y)=\overline{p}_{\nu _{2}}(w;y)$ and (\ref{bi}) easily yields
(\ref{spe}). The latter coincides with the result proved in \cite{BO} and
already recalled in Remark 2.1.

As far as Theorem 2.3 is concerned, by putting $n_{1}=0$, $n_{2}=1,$ the
density of the random time-argument can be expressed as follows:%
\begin{equation}
q_{\nu _{2}}(y,t)=\frac{1}{\lambda t^{\nu _{2}}}\mathcal{W}_{-\nu _{2},1-\nu
_{2}}\left( -\frac{|y|}{\lambda t^{\nu _{2}}}\right) =\overline{v}_{2\nu
_{2}}(y,t),  \label{qu}
\end{equation}%
since in (\ref{bi11}) only the term $r=0$ of the sum survives. To sum up,
the FPP analyzed in \cite{BO} is equal in distribution to the random time
process $N(\mathcal{T}_{\nu }(t))$, whose density can be expressed as a
simple Wright function; on the other hand, in the distributed order case,
the situation is more complicated. The density of the process $N(\widetilde{%
\mathcal{T}}_{\nu _{1},\nu _{2}}(t))$ is written in terms of infinite sums
of Wright functions. Moreover, in the single-order case the density $%
\overline{v}_{2\nu _{2}}(y,t)$ coincides with the folded solution to the
fractional diffusion equation (\ref{due.16}); in the double-order case the
relationship between the density $q_{\nu _{1},\nu _{2}}(y,t)$ and the
fractional diffusion equation of distributed order is more complicated, as
we will prove in the next section.

\

Let us now focus on the probability generating function of the process $%
\widetilde{\mathcal{N}}_{\nu _{1},\nu _{2}}$, which can be expressed in
terms of GML functions (\ref{gml2}), as the following result shows.

\

\noindent \textbf{Theorem 2.4} \ \textit{The probability generating function
}$\widetilde{G}_{\nu _{1},\nu _{2}}(u,t)$\textit{\ of the process }$%
\widetilde{\mathcal{N}}_{\nu _{1},\nu _{2}}$\textit{\ is equal to}%
\begin{align}
\widetilde{G}_{\nu _{1},\nu _{2}}(u,t)& \equiv \sum_{k=0}^{\infty }u^{k}%
\widetilde{p}_{k}^{\nu }(t)  \label{gu} \\
& =\sum_{r=0}^{\infty }\left( -\frac{n_{1}t^{\nu _{2}-\nu _{1}}}{n_{2}}%
\right) ^{r}E_{\nu _{2},(\nu _{2}-\nu _{1})r+1}^{r+1}\left( -\frac{\lambda
(1-u)t^{\nu _{2}}}{n_{2}}\right) +  \notag \\
& -\sum_{r=0}^{\infty }\left( -\frac{n_{1}t^{\nu _{2}-\nu _{1}}}{n_{2}}%
\right) ^{r+1}E_{\nu _{2},(\nu _{2}-\nu _{1})(r+1)+1}^{r+1}\left( -\frac{%
\lambda (1-u)t^{\nu _{2}}}{n_{2}}\right) .  \notag
\end{align}

\noindent \textbf{Proof \ }The Laplace transform of $\widetilde{G}_{\nu
_{1},\nu _{2}}$ can be written, by taking into account\ formula (\ref{bi4}),
as%
\begin{eqnarray}
\mathcal{L}\left\{ \widetilde{G}_{\nu _{1},\nu _{2}}(u,\cdot );\eta \right\}
&=&\sum_{k=0}^{\infty }\frac{u^{k}\lambda ^{k}n_{1}\eta ^{\nu _{1}-1}}{%
(\lambda +n_{1}\eta ^{\nu _{1}}+n_{2}\eta ^{\nu _{2}})^{k+1}}%
+\sum_{k=0}^{\infty }\frac{u^{k}\lambda ^{k}n_{2}\eta ^{\nu _{2}-1}}{%
(\lambda +n_{1}\eta ^{\nu _{1}}+n_{2}\eta ^{\nu _{2}})^{k+1}}  \notag \\
&=&\frac{n_{1}\eta ^{\nu _{1}-1}}{\lambda (1-u)+n_{1}\eta ^{\nu
_{1}}+n_{2}\eta ^{\nu _{2}}}+\frac{n_{2}\eta ^{\nu _{2}-1}}{\lambda
(1-u)+n_{1}\eta ^{\nu _{1}}+n_{2}\eta ^{\nu _{2}}}  \label{gu2} \\
&=&\frac{n_{1}}{n_{2}}\frac{\eta ^{\nu _{1}-1}}{\frac{\lambda }{n_{2}}(1-u)+%
\frac{n_{1}}{n_{2}}\eta ^{\nu _{1}}+\eta ^{\nu _{2}}}+\frac{\eta ^{\nu
_{2}-1}}{\frac{\lambda }{n_{2}}(1-u)+\frac{n_{1}}{n_{2}}\eta ^{\nu
_{1}}+\eta ^{\nu _{2}}}.  \notag
\end{eqnarray}%
By applying formula (26) and (24) of \cite{sa2} to the first \ and second
terms of (\ref{gu2}) respectively and recalling that $\nu _{2}>\nu _{1}$,
the Laplace transform can be inverted as follows:%
\begin{eqnarray*}
\widetilde{G}_{\nu _{1},\nu _{2}}(u,t) &=&\frac{n_{1}}{n_{2}}t^{\nu _{2}-\nu
_{1}}\sum_{r=0}^{\infty }\left( -\frac{n_{1}t^{\nu _{2}-\nu _{1}}}{n_{2}}%
\right) ^{r}E_{\nu _{2},(\nu _{2}-\nu _{1})(r+1)+1}^{r+1}\left( -\frac{%
\lambda (1-u)t^{\nu _{2}}}{n_{2}}\right) + \\
&&+\sum_{r=0}^{\infty }\left( -\frac{n_{1}t^{\nu _{2}-\nu _{1}}}{n_{2}}%
\right) ^{r}E_{\nu _{2},(\nu _{2}-\nu _{1})r+1}^{r+1}\left( -\frac{\lambda
(1-u)t^{\nu _{2}}}{n_{2}}\right) ,
\end{eqnarray*}%
which is equal to (\ref{gu}).\hfill $\square $

\

\noindent \textbf{Remark 2.3 \ }We observe that the infinite sum of GML
functions in (\ref{gu}) coincides with the Fourier transform of the solution
of the diffusion equation (fractional in time and space), analyzed in \cite%
{sa2}, i.e.%
\begin{equation}
\frac{\partial ^{\nu _{2}}}{\partial t^{\nu _{2}}}f(x,t)+\frac{n_{1}}{n_{2}}%
\frac{\partial ^{\nu _{1}}}{\partial t^{\nu _{1}}}f(x,t)=c^{2}\frac{\partial
^{\gamma }}{\partial x^{\gamma }}f(x,t),\quad x\in \mathbb{R},t>0,
\label{gu3}
\end{equation}%
in the special case where $\gamma =0.$ In this case, for $c^{2}=\lambda
(1-u)/n_{2}$, equation (\ref{gu3}) reduces to the time-fractional equation%
\begin{equation}
\frac{\partial ^{\nu _{2}}}{\partial t^{\nu _{2}}}G(u,t)+\frac{n_{1}}{n_{2}}%
\frac{\partial ^{\nu _{1}}}{\partial t^{\nu _{1}}}G(u,t)=\frac{\lambda (1-u)%
}{n_{2}}G(u,t),  \label{gu4}
\end{equation}%
(with initial condition $G(u,0)=1).$ Indeed the probability generating
function $\widetilde{G}_{\nu _{1},\nu _{2}}$ must solve equation (\ref{gu4}%
), as the following steps easily show:%
\begin{eqnarray*}
&&n_{2}\frac{\partial ^{\nu _{2}}}{\partial t^{\nu _{2}}}\widetilde{G}_{\nu
_{1},\nu _{2}}(u,t)+n_{1}\frac{\partial ^{\nu _{1}}}{\partial t^{\nu _{1}}}%
\widetilde{G}_{\nu _{1},\nu _{2}}(u,t) \\
&=&\sum_{k=0}^{\infty }u^{k}\left[ n_{2}\frac{\partial ^{\nu _{2}}}{\partial
t^{\nu _{2}}}\widetilde{p}_{k}^{\nu }(t)+n_{1}\frac{\partial ^{\nu _{1}}}{%
\partial t^{\nu _{1}}}\widetilde{p}_{k}^{\nu }(t)\right] \\
&=&-\lambda \sum_{k=0}^{\infty }u^{k}\widetilde{p}_{k}^{\nu }(t)+\lambda
u\sum_{k=1}^{\infty }u^{k-1}\widetilde{p}_{k-1}^{\nu }(t) \\
&=&-\lambda (1-u)\widetilde{G}_{\nu _{1},\nu _{2}}(u,t).
\end{eqnarray*}

\quad

\noindent \textbf{Remark 2.4 \ }By means of the probability generating
function, we can check that the distribution $\widetilde{p}_{k}^{\nu }(t)$,
sums up to one, for $k=0,1,...$. For $u=1$, formula (\ref{gu}) yields%
\begin{eqnarray}
&&\left. \widetilde{G}_{\nu _{1},\nu _{2}}(u,t)\right\vert
_{u=1}=\sum_{k=0}^{\infty }\widetilde{p}_{k}^{\nu }(t) \\
&=&\sum_{r=0}^{\infty }\left( -\frac{n_{1}t^{\nu _{2}-\nu _{1}}}{n_{2}}%
\right) ^{r}\frac{1}{\Gamma ((\nu _{2}-\nu _{1})r+1)}-\sum_{r=0}^{\infty
}\left( -\frac{n_{1}t^{\nu _{2}-\nu _{1}}}{n_{2}}\right) ^{r+1}\frac{1}{%
\Gamma ((\nu _{2}-\nu _{1})(r+1)+1)}=1,  \notag
\end{eqnarray}%
since only the term $j=0$ in the expression (\ref{gml2}) of the GML function
survives.

Moreover, for $u=0,$ formula (\ref{gu}) gives the probability $\widetilde{p}%
_{0}^{\nu }(t)$ (already obtained in (\ref{pon})):%
\begin{eqnarray}
&&\left. \widetilde{G}_{\nu _{1},\nu _{2}}(u,t)\right\vert _{u=0}=
\label{po} \\
&=&\sum_{r=0}^{\infty }\left( -\frac{n_{1}t^{\nu _{2}-\nu _{1}}}{n_{2}}%
\right) ^{r}E_{\nu _{2},(\nu _{2}-\nu _{1})r+1}^{r+1}\left( -\frac{\lambda
t^{\nu _{2}}}{n_{2}}\right) -\sum_{r=0}^{\infty }\left( -\frac{n_{1}t^{\nu
_{2}-\nu _{1}}}{n_{2}}\right) ^{r+1}E_{\nu _{2},(\nu _{2}-\nu
_{1})(r+1)+1}^{r+1}\left( -\frac{\lambda t^{\nu _{2}}}{n_{2}}\right)  \notag
\\
&=&\widetilde{p}_{0}^{\nu }(t).  \notag
\end{eqnarray}%
We note moreover that, in the special case $n_{1}=0$, $n_{2}=1,$ the
probability generating function reduces to
\begin{equation*}
G_{\nu _{2}}(u,t)=E_{\nu _{2},1}\left( -\lambda (1-u)t^{\nu _{2}}\right)
\end{equation*}%
which coincides with the one obtained for the fractional Poisson process in
\cite{BO}, as expected.

\

We make use of Theorem 2.4 also in the evaluation of the exponential moments
of the process $\widetilde{\mathcal{N}}_{\nu _{1},\nu _{2}}$ and in its
resulting characterization as a Cox process.

\

\noindent \textbf{Theorem 2.5} \ \textit{The factorial moments of the
process }$\widetilde{\mathcal{N}}_{\nu _{1},\nu _{2}},$\textit{\ with
distribution }$\widetilde{p}_{k}^{\nu }(t)$\textit{\ and probability
generating function }$\widetilde{G}_{\nu _{1},\nu _{2}}(u,t)$\textit{\ given
in (\ref{gu}), are equal to}%
\begin{eqnarray}
&&\mathbb{E}\left[ \widetilde{\mathcal{N}}_{\nu _{1},\nu _{2}}(t)\left(
\widetilde{\mathcal{N}}_{\nu _{1},\nu _{2}}(t)-1\right) ...(\widetilde{%
\mathcal{N}}_{\nu _{1},\nu _{2}}(t)-k+1)\right]  \label{fac} \\
&=&\frac{n_{1}\lambda ^{k}t^{\nu _{2}k+\nu _{2}-\nu _{1}}}{n_{2}^{k+1}}%
k!E_{\nu _{2}-\nu _{1},\nu _{2}k+\nu _{2}-\nu _{1}+1}^{k+1}\left( -\frac{%
n_{1}t^{\nu _{2}-\nu _{1}}}{n_{2}}\right) +\frac{\lambda ^{k}t^{\nu _{2}k}}{%
n_{2}^{k}}k!E_{\nu _{2}-\nu _{1},\nu _{2}k+1}^{k+1}\left( -\frac{n_{1}t^{\nu
_{2}-\nu _{1}}}{n_{2}}\right) .  \notag
\end{eqnarray}%
\textit{Moreover }$\widetilde{\mathcal{N}}_{\nu _{1},\nu _{2}}$\textit{\ is
a Cox process with directing measure }$\Lambda \left( \left( 0,t\right]
\right) \equiv \widetilde{\mathcal{T}}_{\nu _{1},\nu _{2}}(t),$\textit{\
endowed with density }$q_{\nu _{1},\nu _{2}}(y,t).$

\noindent \textbf{Proof \ }We take the $k$-th derivatives of $\widetilde{G}%
_{\nu _{1},\nu _{2}}$ with respect to $u$:%
\begin{eqnarray*}
&&\frac{\partial ^{k}}{\partial u^{k}}\widetilde{G}_{\nu _{1},\nu _{2}}(u,t)
\\
&=&\sum_{r=0}^{\infty }\left( -\frac{n_{1}t^{\nu _{2}-\nu _{1}}}{n_{2}}%
\right) ^{r}\frac{1}{r!}\sum_{j=0}^{\infty }\frac{(r+j)!}{j!\Gamma (\nu
_{2}j+(\nu _{2}-\nu _{1})r+1)}\frac{\left( -\frac{\lambda t^{\nu _{2}}}{n_{2}%
}\right) ^{j}(-1)^{k}(1-u)^{j-k}j!}{(j-k)!}+ \\
&&-\sum_{r=0}^{\infty }\left( -\frac{n_{1}t^{\nu _{2}-\nu _{1}}}{n_{2}}%
\right) ^{r+1}\frac{1}{r!}\sum_{j=0}^{\infty }\frac{(r+j)!}{j!\Gamma (\nu
_{2}j+(\nu _{2}-\nu _{1})(r+1)+1)}\frac{\left( -\frac{\lambda t^{\nu _{2}}}{%
n_{2}}\right) ^{j}(-1)^{k}(1-u)^{j-k}j!}{(j-k)!},
\end{eqnarray*}%
which, for $u=1$, becomes%
\begin{eqnarray}
\left. \frac{\partial ^{k}}{\partial u^{k}}\widetilde{G}_{\nu _{1},\nu
_{2}}(u,t)\right\vert _{u=1} &=&\sum_{r=0}^{\infty }\left( -\frac{%
n_{1}t^{\nu _{2}-\nu _{1}}}{n_{2}}\right) ^{r}\frac{1}{r!}\frac{(r+k)!\left(
\frac{\lambda t^{\nu _{2}}}{n_{2}}\right) ^{k}}{\Gamma (\nu _{2}k+(\nu
_{2}-\nu _{1})r+1)}  \label{ra} \\
&&-\sum_{r=0}^{\infty }\left( -\frac{n_{1}t^{\nu _{2}-\nu _{1}}}{n_{2}}%
\right) ^{r+1}\frac{1}{r!}\frac{(r+k)!\left( \frac{\lambda t^{\nu _{2}}}{%
n_{2}}\right) ^{k}}{\Gamma (\nu _{2}k+(\nu _{2}-\nu _{1})(r+1)+1)}.  \notag
\end{eqnarray}%
Formula (\ref{ra}) can be written as (\ref{fac}), by multiplying and
dividing for $k!.$

In order to prove that $\widetilde{\mathcal{N}}_{\nu _{1},\nu _{2}}(t),t>0$
is a Cox process with directing measure equal to $\Lambda \left( \left( 0,t%
\right] \right) \equiv \widetilde{\mathcal{T}}_{\nu _{1},\nu _{2}}(t)$, we
adopt the characterization of Cox processes by its factorial moments. Indeed
it is proved in \cite{La} that for a Cox process they must coincide with the
ordinary moments of its directing measure. Our goal is to show that this
equivalence holds for $\widetilde{\mathcal{N}}_{\nu _{1},\nu _{2}}$ and for
the density $q_{\nu _{1},\nu _{2}}(y,t)$ of its time argument, i.e. that%
\begin{equation}
\mathbb{E}\left[ \widetilde{\mathcal{T}}_{\nu _{1},\nu _{2}}(t)\right]
^{k}=\int_{0}^{+\infty }y^{k}q_{\nu _{1},\nu _{2}}(y,t)dy  \label{e1}
\end{equation}%
coincides with (\ref{fac}). We start by taking the Laplace transform of (\ref%
{bi11}), which reads%
\begin{eqnarray}
\mathcal{L}\left\{ q_{\nu _{1}.\nu _{2}}(y,\cdot );\eta \right\} &=&\frac{%
n_{1}}{\lambda }\sum_{r=0}^{\infty }\frac{\left( -n_{2}|y|/\lambda \right)
^{r}}{r!}\sum_{l=0}^{\infty }\frac{\left( -n_{1}|y|/\lambda \right) ^{l}}{%
l!\Gamma (-\nu _{1}l+1-\nu _{2}r-\nu _{1})}\frac{\Gamma (1-\nu _{1}l-\nu
_{2}r-\nu _{1})}{\eta ^{1-\nu _{1}l-\nu _{2}r-\nu _{1}}}+  \notag \\
&&+\frac{n_{2}}{\lambda }\sum_{r=0}^{\infty }\frac{\left( -n_{1}|y|/\lambda
\right) ^{r}}{r!}\sum_{l=0}^{\infty }\frac{\left( -n_{2}|y|/\lambda \right)
^{l}}{l!\Gamma (-\nu _{2}l+1-\nu _{1}r-\nu _{2})}\frac{\Gamma (1-\nu
_{2}l-\nu _{1}r-\nu _{2})}{\eta ^{1-\nu _{2}l-\nu _{1}r-\nu _{2}}}  \notag \\
&=&\frac{n_{1}e^{-n_{1}\eta ^{\nu _{1}}|y|/\lambda }}{\lambda \eta ^{1-\nu
_{1}}}\sum_{r=0}^{\infty }\frac{\left( -n_{2}\eta ^{\nu _{2}}|y|/\lambda
\right) ^{r}}{r!}+\frac{n_{2}e^{-n_{2}\eta ^{\nu _{2}}|y|/\lambda }}{\lambda
\eta ^{1-\nu _{2}}}\sum_{r=0}^{\infty }\frac{\left( -n_{1}\eta ^{\nu
_{1}}|y|/\lambda \right) ^{r}}{r!}  \notag \\
&=&\frac{n_{1}e^{-\left( n_{1}\eta ^{\nu _{1}}+n_{2}\eta ^{\nu _{2}}\right)
|y|/\lambda }}{\lambda \eta ^{1-\nu _{1}}}+\frac{n_{2}e^{-(n_{2}\eta ^{\nu
_{2}}+n_{1}\eta ^{\nu _{1}})|y|/\lambda }}{\lambda \eta ^{1-\nu _{2}}},
\label{e2}
\end{eqnarray}%
so that the Laplace transform of (\ref{e1}) becomes%
\begin{eqnarray}
&&\mathcal{L}\left\{ \mathbb{E}\left[ \widetilde{\mathcal{T}}_{\nu _{1},\nu
_{2}}(\cdot )\right] ^{k};\eta \right\}  \label{fac3} \\
&=&\frac{n_{1}}{\lambda \eta ^{1-\nu _{1}}}\int_{0}^{+\infty
}y^{k}e^{-\left( n_{1}\eta ^{\nu _{1}}+n_{2}\eta ^{\nu _{2}}\right)
y/\lambda }dy+\frac{n_{2}}{\lambda \eta ^{1-\nu _{2}}}\int_{0}^{+\infty
}y^{k}e^{-\left( n_{1}\eta ^{\nu _{1}}+n_{2}\eta ^{\nu _{2}}\right)
y/\lambda }dy  \notag \\
&=&\frac{n_{1}\lambda ^{k}\eta ^{\nu _{1}-1}k!}{\left( n_{1}\eta ^{\nu
_{1}}+n_{2}\eta ^{\nu _{2}}\right) ^{k+1}}+\frac{n_{2}\lambda ^{k}\eta ^{\nu
_{2}-1}k!}{\left( n_{1}\eta ^{\nu _{1}}+n_{2}\eta ^{\nu _{2}}\right) ^{k+1}}.
\notag
\end{eqnarray}%
We now take the Laplace transform of (\ref{fac}), by applying (\ref{pra}):%
\begin{eqnarray*}
&&\mathcal{L}\left\{ \mathbb{E}\left[ \widetilde{\mathcal{N}}_{\nu _{1},\nu
_{2}}(\cdot )...(\widetilde{\mathcal{N}}_{\nu _{1},\nu _{2}}(\cdot )-k+1)%
\right] ;\eta \right\} \\
&=&\frac{n_{1}\lambda ^{k}}{n_{2}^{k+1}}k!\frac{\eta ^{-\nu _{1}k-1}}{\left(
\eta ^{\nu _{2}-\nu _{1}}+\frac{n_{1}}{n_{2}}\right) ^{k+1}}+\frac{\lambda
^{k}}{n_{2}^{k}}k!\frac{\eta ^{\nu _{2}-\nu _{1}-\nu _{1}k-1}}{\left( \eta
^{\nu _{2}-\nu _{1}}+\frac{n_{1}}{n_{2}}\right) ^{k+1}}.
\end{eqnarray*}%
It is simply verified that the last expression coincides with (\ref{fac3}).%
\hfil$\square $

\

\noindent \textbf{Remark 2.5 \ }For $k=1$ we get from (\ref{fac}) the
expected value of $\widetilde{\mathcal{N}}_{\nu _{1},\nu _{2}}$:%
\begin{eqnarray}
&&\mathbb{E}\widetilde{\mathcal{N}}_{\nu _{1},\nu _{2}}(t)  \label{exp} \\
&=&\frac{n_{1}\lambda t^{2\nu _{2}-\nu _{1}}}{n_{2}^{2}}E_{\nu _{2}-\nu
_{1},2\nu _{2}-\nu _{1}+1}^{2}\left( -\frac{n_{1}t^{\nu _{2}-\nu _{1}}}{n_{2}%
}\right) +\frac{\lambda t^{\nu _{2}}}{n_{2}}E_{\nu _{2}-\nu _{1},\nu
_{2}+1}^{2}\left( -\frac{n_{1}t^{\nu _{2}-\nu _{1}}}{n_{2}}\right)  \notag \\
&=&\frac{n_{1}\lambda t^{2\nu _{2}-\nu _{1}}}{n_{2}^{2}}\sum_{j=0}^{\infty
}\left( -\frac{n_{1}t^{\nu _{2}-\nu _{1}}}{n_{2}}\right) ^{j}\frac{j+1}{%
\Gamma (\nu _{2}j-\nu _{1}j+2\nu _{2}-\nu _{1}+1)}+  \notag \\
&&+\frac{\lambda t^{\nu _{2}}}{n_{2}}\sum_{j=0}^{\infty }\left( -\frac{%
n_{1}t^{\nu _{2}-\nu _{1}}}{n_{2}}\right) ^{j}\frac{j+1}{\Gamma (\nu
_{2}j-\nu _{1}j+\nu _{2}+1)}  \notag \\
&=&-\frac{\lambda t^{\nu _{2}}}{n_{2}}\sum_{l=1}^{\infty }\left( -\frac{%
n_{1}t^{\nu _{2}-\nu _{1}}}{n_{2}}\right) ^{l}\frac{l}{\Gamma ((\nu _{2}-\nu
_{1})l+\nu _{2}+1)}+  \notag \\
&&+\frac{\lambda t^{\nu _{2}}}{n_{2}}\sum_{j=0}^{\infty }\left( -\frac{%
n_{1}t^{\nu _{2}-\nu _{1}}}{n_{2}}\right) ^{j}\frac{j+1}{\Gamma ((\nu
_{2}-\nu _{1})j+\nu _{2}+1)}  \notag \\
&=&\frac{\lambda t^{\nu _{2}}}{n_{2}}\sum_{j=0}^{\infty }\left( -\frac{%
n_{1}t^{\nu _{2}-\nu _{1}}}{n_{2}}\right) ^{j}\frac{1}{\Gamma ((\nu _{2}-\nu
_{1})j+\nu _{2}+1)}  \notag \\
&=&\frac{\lambda t^{\nu _{2}}}{n_{2}}E_{\nu _{2}-\nu _{1},\nu _{2}+1}\left( -%
\frac{n_{1}t^{\nu _{2}-\nu _{1}}}{n_{2}}\right) .  \notag
\end{eqnarray}%
Now consider again the particular case $n_{1}=0$, $n_{2}=1$; formula (\ref%
{fac}) reduces, in this case, to
\begin{equation}
\mathbb{E}\left[ \mathcal{N}_{\nu _{2}}(t)\left( \mathcal{N}_{\nu
_{2}}(t)-1\right) ...(\mathcal{N}_{\nu _{2}}(t)-k+1)\right] =\frac{\lambda
^{k}t^{\nu _{2}k}k!}{\Gamma \left( \nu _{2}k+1\right) },  \notag
\end{equation}%
which coincides with the factorial moments of the FPP obtained in \cite{BO}.
Analogously the expected value given in (\ref{exp}) reduces (for $n_{1}=0$, $%
n_{2}=1$) to
\begin{equation}
\mathbb{E}\mathcal{N}_{\nu _{2}}(t)=\frac{\lambda t^{\nu _{2}}}{\Gamma
\left( \nu _{2}+1\right) },  \label{exp2}
\end{equation}%
as expected. We observe that, in the distributed order case analyzed here,
both the factorial moments and the expected value of $\widetilde{\mathcal{N}}%
_{\nu _{1},\nu _{2}}$ are expressed in terms of Mittag-Leffler functions;
for $k=1$ it is a two-parameter Mittag-Leffler function, while, for $k>1,$
we need a GML function with third parameter equal to $k+1.$ This is
analogously true, in view of Theorem 2.5, for the $k$-th order moments of
the time argument $\widetilde{\mathcal{T}}_{\nu _{1},\nu _{2}}(t).$

\

We concentrate now on the renewal property of $\widetilde{\mathcal{N}}_{\nu
_{1},\nu _{2}}$: more precisely, we obtain the density of the waiting-time
of the $k$-th event $\widetilde{f}_{k}^{\nu }(t)$ (or, more exactly, its
Laplace transform) and that of the interarrival times $\widetilde{f}%
_{1}^{\nu }(t)$. The latter is expressed again by means of infinite sums of
GML functions. The same is true for the survival probability $\widetilde{%
\Psi }^{\nu }(t)$. We remark that $\widetilde{f}_{k}^{\nu }(t)$ can be
expressed as the $k$-th convolution of $\widetilde{f}_{1}^{\nu }(t)$ and
this implies that the process $\widetilde{\mathcal{N}}_{\nu _{1},\nu _{2}}$
is a renewal, since the waiting-time of the $k$-th event $T_{k}\equiv \inf
\left\{ t>0:\widetilde{\mathcal{N}}_{\nu _{1},\nu _{2}}(t)=k\right\} $ is
given by the sum of $k$ independent and identically distributed interarrival
times $U_{j},j=1,...k$.

\

\noindent \textbf{Theorem 2.6 \ }\textit{The Laplace transform of the
density }$\widetilde{f}_{k}^{\nu }(t)=\Pr \left\{ T_{k}\in dt\right\} $%
\textit{\ of the }$k$\textit{-th event waiting-time }$T_{k}$\textit{, is
equal to}%
\begin{equation}
\mathcal{L}\left\{ \widetilde{f}_{k}^{\nu };\eta \right\} =\left( \frac{%
\lambda }{\lambda +n_{1}\eta ^{\nu _{1}}+n_{2}\eta ^{\nu _{2}}}\right)
^{k},\quad k\geq 1.  \label{tem1}
\end{equation}%
\textit{The density of the interarrival time }$U_{j}$\textit{\ is equal to }$%
\widetilde{f}_{1}^{\nu },$\textit{\ for any }$j=1,2,...$\textit{and can be
written as}%
\begin{equation}
\widetilde{f}_{1}^{\nu }(t)=\frac{\lambda }{n_{2}}t^{\nu
_{2}-1}\sum_{r=0}^{\infty }\left( -\frac{n_{1}t^{\nu _{2}-\nu _{1}}}{n_{2}}%
\right) ^{r}E_{\nu _{2},\nu _{2}+(\nu _{2}-\nu _{1})r}^{r+1}\left( -\frac{%
\lambda t^{\nu _{2}}}{n_{2}}\right) .  \label{tem2}
\end{equation}%
\textit{Alternatively }%
\begin{equation}
\widetilde{f}_{1}^{\nu }(t)=\int_{0}^{+\infty }f_{1}(s)g_{\nu _{1},\nu
_{2}}(s,t)ds=\int_{0}^{+\infty }e^{-s}g_{\nu _{1},\nu _{2}}(s,t)ds,
\label{tem3}
\end{equation}%
\textit{where }$f_{1}$\textit{\ denotes the interarrival-time density of the
Poisson process }$N$\textit{\ (i.e. }$f_{1}(t)=e^{-t})$\textit{\ and }$%
g_{\nu _{1},\nu _{2}}$\textit{\ is given in (\ref{bi9}). Then }$\widetilde{%
\mathcal{N}}_{\nu _{1},\nu _{2}}$\textit{\ is a renewal process with renewal
function given by}%
\begin{equation}
\widetilde{m}^{\nu }(t)=\frac{\lambda t^{\nu _{2}}}{n_{2}}E_{\nu _{2}-\nu
_{1},\nu _{2}+1}\left( -\frac{n_{1}t^{\nu _{2}-\nu _{1}}}{n_{2}}\right) .
\label{emme}
\end{equation}%
\textit{The survival probability }$\widetilde{\Psi }^{\nu }(t)\equiv \Pr
\left\{ U_{1}>t\right\} $\textit{\ can be expressed as}%
\begin{eqnarray}
&&\widetilde{\Psi }^{\nu }(t)=1-\int_{0}^{t}\widetilde{f}_{1}^{\nu }(s)ds
\label{tem4} \\
&=&\sum_{r=0}^{\infty }\left( -\frac{n_{1}t^{\nu _{2}-\nu _{1}}}{n_{2}}%
\right) ^{r}E_{\nu _{2},(\nu _{2}-\nu _{1})r+1}^{r+1}\left( -\frac{\lambda
t^{\nu _{2}}}{n_{2}}\right) -\sum_{r=0}^{\infty }\left( -\frac{n_{1}t^{\nu
_{2}-\nu _{1}}}{n_{2}}\right) ^{r+1}E_{\nu _{2},(\nu _{2}-\nu
_{1})(r+1)+1}^{r+1}\left( -\frac{\lambda t^{\nu _{2}}}{n_{2}}\right) ,
\notag
\end{eqnarray}%
\textit{which solves the relaxation equation of distributed order}%
\begin{equation}
n_{1}\frac{\partial ^{\nu _{1}}}{\partial t^{\nu _{1}}}\Psi (t)+n_{2}\frac{%
\partial ^{\nu _{2}}}{\partial t^{\nu _{2}}}\Psi (t)=-\lambda \Psi (t),
\label{tem7}
\end{equation}%
\textit{(with initial condition }$\Psi (0)=1).$

\noindent \textbf{Proof \ }Formula (\ref{tem1}) easily follows from the
following relationship%
\begin{eqnarray*}
&&\mathcal{L}\left\{ \widetilde{p}_{k}^{\nu };\eta \right\}
=\int_{0}^{+\infty }e^{-\eta t}\Pr \left\{ \widetilde{\mathcal{N}}_{\nu
_{1},\nu _{2}}(t)=k\right\} dt \\
&=&\int_{0}^{+\infty }e^{-\eta t}\left[ \Pr \left\{ T_{k}<t\right\} -\Pr
\left\{ T_{k+1}<t\right\} \right] dt \\
&=&\frac{1}{\eta }\left[ \mathcal{L}\left\{ \widetilde{f}_{k}^{\nu };\eta
\right\} -\mathcal{L}\left\{ \widetilde{f}_{k+1}^{\nu };\eta \right\} \right]
,
\end{eqnarray*}%
used together with (\ref{bi4}). The Laplace transform of the density of the
first interarrival time $U_{1}$ is equal to (\ref{tem1}) for $k=1$:%
\begin{equation}
\mathcal{L}\left\{ \widetilde{f}_{1}^{\nu };\eta \right\} =\frac{\lambda }{%
\lambda +n_{1}\eta ^{\nu _{1}}+n_{2}\eta ^{\nu _{2}}},  \label{tem5}
\end{equation}%
and thus the density of the $k$-th event waiting time $\widetilde{f}%
_{k}^{\nu }$ is expressed as the\ $k$-fold convolution of $\widetilde{f}%
_{1}^{\nu }.$ This proves that $\widetilde{\mathcal{N}}_{\nu _{1},\nu _{2}}$
is a renewal process; its renewal function has been already calculated in (%
\ref{exp}). As a check we show that the well-known relationship between the
Laplace transforms of $\widetilde{m}^{\nu }(t)$ and $\widetilde{f}_{1}^{\nu
} $ holds in this case:%
\begin{eqnarray*}
\mathcal{L}\left\{ \widetilde{m}^{\nu };\eta \right\} &=&\frac{\lambda }{%
\eta \left[ n_{1}\eta ^{\nu _{1}}+n_{2}\eta ^{\nu _{2}}\right] } \\
&=&\frac{\frac{\lambda }{\lambda +n_{1}\eta ^{\nu _{1}}+n_{2}\eta ^{\nu _{2}}%
}}{\eta \left[ 1-\frac{\lambda }{\lambda +n_{1}\eta ^{\nu _{1}}+n_{2}\eta
^{\nu _{2}}}\right] }=\frac{\mathcal{L}\left\{ \widetilde{f}_{1}^{\nu };\eta
\right\} }{\eta \left[ 1-\mathcal{L}\left\{ \widetilde{f}_{1}^{\nu };\eta
\right\} \right] }.
\end{eqnarray*}

The Laplace transform (\ref{tem5}) can be inverted by applying formula (27)
of \cite{sa2}, for $\alpha =\nu _{2},$ $a=n_{1}/n_{2},$ $\beta =\nu _{1}$
and $b=\lambda /n_{2},$ thus giving (\ref{tem2}). We can rewrite moreover (%
\ref{tem5}) as follows:%
\begin{eqnarray}
\mathcal{L}\left\{ \widetilde{f}_{1}^{\nu };\eta \right\} &=&\frac{\lambda }{%
\lambda +n_{1}\eta ^{\nu _{1}}+n_{2}\eta ^{\nu _{2}}} \\
&=&\mathcal{L}\left\{ f_{1};\frac{n_{1}}{\lambda }\eta ^{\nu _{1}}+\frac{%
n_{2}}{\lambda }\eta ^{\nu _{2}}\right\} ,  \notag
\end{eqnarray}%
where $f_{1}(t),t>0,$ is again the density of the interarrival times for the
Poisson process; hence%
\begin{equation}
\mathcal{L}\left\{ \widetilde{f}_{1}^{\nu };\eta \right\} =\int_{0}^{+\infty
}f_{1}(t)e^{-(\frac{n_{1}}{\lambda }\eta ^{\nu _{1}}+\frac{n_{2}}{\lambda }%
\eta ^{\nu _{2}})t}dt.  \label{tem6}
\end{equation}%
By inverting the Laplace transform (\ref{tem6}), taking into account (\ref%
{bi}) and (\ref{bi9}), we get (\ref{tem3}).

We take the derivative of (\ref{tem4}) and we show that $-\frac{d}{dt}%
\widetilde{\Psi }^{\nu }(t)=\widetilde{f}_{1}^{\nu }(t)$ given in (\ref{tem2}%
): indeed%
\begin{eqnarray*}
&&-\frac{d}{dt}\widetilde{\Psi }^{\nu }(t) \\
&=&-(\nu _{2}-\nu _{1})\sum_{r=0}^{\infty }r\left( -\frac{n_{1}}{n_{2}}%
\right) ^{r}t^{(\nu _{2}-\nu _{1})r-1}E_{\nu _{2},(\nu _{2}-\nu
_{1})r+1}^{r+1}\left( -\frac{\lambda t^{\nu _{2}}}{n_{2}}\right) + \\
&&-\sum_{r=0}^{\infty }\left( -\frac{n_{1}t^{\nu _{2}-\nu _{1}}}{n_{2}}%
\right) ^{r}\frac{1}{r!}\sum_{j=0}^{\infty }\frac{(r+j)!\left( -\frac{%
\lambda }{n_{2}}\right) ^{j}j\nu _{2}t^{\nu _{2}j-1}}{j!\Gamma (\nu
_{2}j+(\nu _{2}-\nu _{1})r+1)}+ \\
&&+(\nu _{2}-\nu _{1})\sum_{r=0}^{\infty }(r+1)\left( -\frac{n_{1}}{n_{2}}%
\right) ^{r+1}t^{(\nu _{2}-\nu _{1})(r+1)-1}E_{\nu _{2},(\nu _{2}-\nu
_{1})(r+1)+1}^{r+1}\left( -\frac{\lambda t^{\nu _{2}}}{n_{2}}\right) + \\
&&+\sum_{r=0}^{\infty }\left( -\frac{n_{1}t^{\nu _{2}-\nu _{1}}}{n_{2}}%
\right) ^{r+1}\frac{1}{r!}\sum_{j=0}^{\infty }\frac{(r+j)!\left( -\frac{%
\lambda }{n_{2}}\right) ^{j}j\nu _{2}t^{\nu _{2}j-1}}{j!\Gamma (\nu
_{2}j+(\nu _{2}-\nu _{1})(r+1)+1)} \\
&=&-\frac{1}{t}\sum_{r=0}^{\infty }\left( -\frac{n_{1}t^{\nu _{2}-\nu _{1}}}{%
n_{2}}\right) ^{r}E_{\nu _{2},(\nu _{2}-\nu _{1})r}^{r+1}\left( -\frac{%
\lambda t^{\nu _{2}}}{n_{2}}\right) +\frac{1}{t}\sum_{r=0}^{\infty }\left( -%
\frac{n_{1}t^{\nu _{2}-\nu _{1}}}{n_{2}}\right) ^{r+1}E_{\nu _{2},(\nu
_{2}-\nu _{1})(r+1)}^{r+1}\left( -\frac{\lambda t^{\nu _{2}}}{n_{2}}\right) .
\end{eqnarray*}%
The latter expression can be shown to coincide with (\ref{tem2}).

By noting that (\ref{tem4}) is equal to (\ref{gu}) for $u=0$, it is
immediately proved that $\widetilde{\Psi }^{\nu }$ solves equation (\ref{gu4}%
) for $u=0,$ i.e. equation (\ref{tem7}). Alternatively, it is easy to check
that the Laplace transform of (\ref{tem4}) is given by%
\begin{equation*}
\mathcal{L}\left\{ \widetilde{\Psi }^{\nu };\eta \right\} =\frac{n_{1}\eta
^{\nu _{1}-1}+n_{2}\eta ^{\nu _{2}-1}}{\lambda +n_{1}\eta ^{\nu
_{1}}+n_{2}\eta ^{\nu _{2}}},
\end{equation*}%
which coincides with the solution to the Laplace transform of equation (\ref%
{tem7}).

\hfill $\square $

\

\noindent \textbf{Remark 2.6 \ }In the special case $n_{1}=0$, $n_{2}=1$,
from (\ref{tem2}) we retrieve the density of the interarrival times of the
fractional Poisson process (see \cite{BO}):%
\begin{equation}
f_{1}^{\nu }(t)=\lambda t^{\nu _{2}-1}E_{\nu _{2},\nu _{2}}\left( -\lambda
t^{\nu _{2}}\right) .  \label{rai9}
\end{equation}%
Likewise the survival probability (\ref{tem4}) reduces to%
\begin{equation}
\Psi ^{\nu }(t)=E_{\nu _{2},1}\left( -\lambda t^{\nu _{2}}\right) .  \notag
\end{equation}

\

It is interesting to analyze the asymptotic behavior of the waiting time
densities and of the renewal function and to compare these expressions with
the corresponding formulas obtained for the fractional Poisson process. To
this purpose we need to prove the following integral representation for the
GML function:%
\begin{equation}
E_{\nu ,\beta }^{k}(-ct^{\nu })=\frac{t^{1-\beta }}{2\pi i}\int_{0}^{+\infty
}r^{\nu k-\beta }e^{-rt}\frac{e^{i\pi \beta }(r^{\nu }+ce^{-i\pi \nu
})^{k}-e^{-i\pi \beta }(r^{\nu }+ce^{i\pi \nu })^{k}}{\left[ r^{2\nu
}+2r^{\nu }c\cos (\nu \pi )+c^{2}\right] ^{k}}dr.  \label{rai8}
\end{equation}%
We start by checking that, for $k=1,$ formula (\ref{rai8}) coincides with
the form given for $E_{\nu ,\beta }(-t^{\nu })$ in \cite{BO}, i.e.%
\begin{equation}
E_{\nu ,\beta }(-ct^{\nu })=\frac{t^{1-\beta }}{\pi }\int_{0}^{+\infty
}r^{\nu -\beta }e^{-rt}\frac{r^{\nu }\sin (\beta \pi )+c\sin ((\beta -\nu
)\pi )}{r^{2\nu }+2r^{\nu }c\cos (\nu \pi )+c^{2}}dr.  \label{rai10}
\end{equation}%
In order to prove formula (\ref{rai8}) we multiply and divide the $m$-th
term in the series expression of $E_{\nu ,\beta }^{k}(-ct^{\nu })$ for $\sin
((\beta +\nu m)\pi )/\pi $ and apply again the reflection formula of the
Gamma function, as follows%
\begin{eqnarray}
&&E_{\nu ,\beta }^{k}(-ct^{\nu })  \label{asy} \\
&=&\frac{1}{(k-1)!}\sum_{m=0}^{\infty }\frac{(m+k-1)!(-ct^{\nu })^{m}}{%
m!\Gamma \left( \nu m+\beta \right) }\frac{\sin ((\beta +\nu m)\pi )}{\pi }%
\Gamma (1-\nu m-\beta )\Gamma (\nu m+\beta )  \notag \\
&=&\frac{t^{1-\beta }}{\pi (k-1)!}\sum_{m=0}^{\infty }\frac{(m+k-1)!(-c)^{m}%
}{m!\Gamma \left( \nu m+\beta \right) }\sin ((\beta +\nu m)\pi
)\int_{0}^{\infty }e^{-rt}r^{-\nu m-\beta }dr\int_{0}^{\infty }e^{-y}y^{\nu
m+\beta -1}dy  \notag \\
&=&\frac{t^{1-\beta }}{\pi (k-1)!}\int_{0}^{\infty
}e^{-rt}dr\int_{0}^{\infty }e^{-ry}y^{\beta -1}\sum_{m=0}^{\infty }\frac{%
(m+k-1)!(-cy^{\nu })^{m}}{m!\Gamma \left( \nu m+\beta \right) }\frac{e^{i\pi
\nu m+i\pi \beta }-e^{-i\pi \nu m-i\pi \beta }}{2i}dy  \notag \\
&=&\frac{t^{1-\beta }}{2\pi i(k-1)!}\int_{0}^{\infty }dy\int_{0}^{\infty
}e^{-r(y+t)}y^{\beta -1}\left[ e^{i\pi \beta }E_{\nu ,\beta }^{k}(-cy^{\nu
}e^{i\pi \nu })-e^{-i\pi \beta }E_{\nu ,\beta }^{k}(-cy^{\nu }e^{-i\pi \nu })%
\right] dr  \notag \\
&=&\left[ \text{by (\ref{pra})}\right]  \notag \\
&=&\frac{t^{1-\beta }}{2\pi i}\int_{0}^{\infty }e^{-rt}\left[ e^{i\pi \beta }%
\frac{r^{\nu k-\beta }}{(r^{\nu }+ce^{i\pi \nu })^{k}}-e^{-i\pi \beta }\frac{%
r^{\nu k-\beta }}{(r^{\nu }+ce^{-i\pi \nu })^{k}}\right] dr.  \notag
\end{eqnarray}%
This coincides with (\ref{rai8}). The asymptotic behavior of $E_{\nu ,\beta
}^{k}(-ct^{\nu })$ can be obtained from (\ref{rai8}) and reads, for $%
t\rightarrow \infty $:%
\begin{eqnarray}
E_{\nu ,\beta }^{k}(-ct^{\nu }) &=&\frac{t^{-\nu k}}{2\pi i}\int_{0}^{\infty
}e^{-z}z^{\nu k-\beta }\left[ \frac{e^{i\pi \beta }(\frac{z^{\nu }}{t^{\nu }}%
+ce^{-i\pi \nu })^{k}-e^{-i\pi \beta }(\frac{z^{\nu }}{t^{\nu }}+ce^{i\pi
\nu })^{k}}{(\frac{z^{2\nu }}{t^{2\nu }}+2c\frac{z^{\nu }}{t^{\nu }}\cos \pi
\nu +c^{2})^{k}}\right] dz  \notag \\
&=&\frac{\Gamma (1-\beta +\nu k)}{\pi c^{k}t^{\nu k}}\sin (\pi (\beta -\nu
k))+o(t^{-\nu k})  \label{asy3} \\
&=&\frac{1}{c^{k}t^{\nu k}\Gamma (\beta -\nu k)}+o(t^{-\nu k}),\qquad
t\rightarrow \infty .  \notag
\end{eqnarray}%
For $k=1$, formula (\ref{asy3}) reduces to the one holding for the
Mittag-Leffler function, which can be deduced by (\ref{rai10}), i.e.%
\begin{equation}
E_{\nu ,\beta }(-ct^{\nu })=\frac{1}{ct^{\nu }\Gamma (\beta -\nu )}%
+o(t^{-\nu }),\qquad t\rightarrow \infty .  \label{asy4}
\end{equation}%
For $t\rightarrow 0$, we get instead that%
\begin{eqnarray}
\lim_{t\rightarrow 0^{+}}E_{\nu ,\beta }^{k}(-ct^{\nu })
&=&\lim_{t\rightarrow 0^{+}}\frac{1}{2\pi i}\int_{0}^{\infty }e^{-z}z^{\nu
k-\beta }\left[ \frac{e^{i\pi \beta }(z^{\nu }+ce^{-i\pi \nu }t^{\nu
})^{k}-e^{-i\pi \beta }(z^{\nu }+ce^{i\pi \nu }t^{\nu })^{k}}{(z^{2\nu
}+2cz^{\nu }t^{\nu }\cos \pi \nu +c^{2}t^{2\nu })^{k}}\right] dz  \notag \\
&=&\frac{1}{\pi }\int_{0}^{\infty }e^{-z}z^{2\nu k-\beta }\frac{\sin (\pi
\beta )}{z^{2\nu k}}dz=\frac{\Gamma (1-\beta )}{\pi }\sin (\pi \beta )=\frac{%
1}{\Gamma (\beta )}.  \label{asy2}
\end{eqnarray}%
The asymptotic behavior for small $t$ can be deduced directly by the series
expression of $E_{\nu ,\beta }^{k}(-ct^{\nu })$: indeed we get, for $0\leq
t<<1,$%
\begin{equation*}
E_{\nu ,\beta }^{k}(-ct^{\nu })\simeq \frac{1}{\Gamma (\beta )}-\frac{%
ct^{\nu }k}{\Gamma (\beta +\nu )},
\end{equation*}%
which reduces, for $k=1,$ to the well-known expression (see \cite{marab},
formula (3.13)).

The interarrival-time density (\ref{tem2}) can be rewritten, by applying (%
\ref{asy}), as

\begin{eqnarray}
&&\widetilde{f}_{1}^{\nu }(t)  \label{inter} \\
&=&\frac{\lambda t^{\nu _{2}-1}}{n_{2}}\sum_{r=0}^{\infty }\left( -\frac{%
n_{1}t^{\nu _{2}-\nu _{1}}}{n_{2}}\right) ^{r}\frac{t^{1-\nu _{2}-(\nu
_{2}-\nu _{1})r}}{2\pi i}\cdot  \notag \\
&&\cdot \int_{0}^{\infty }e^{-zt}z^{\nu _{2}(r+1)-\nu _{2}-(\nu _{2}-\nu
_{1})r}\left[ \frac{e^{i\pi \nu _{2}+i\pi (\nu _{2}-\nu _{1})r}}{(z^{\nu
_{2}}+\frac{\lambda }{n_{2}}e^{i\pi \nu _{2}})^{r+1}}-\frac{e^{-i\pi \nu
_{2}-i\pi (\nu _{2}-\nu _{1})r}}{(z^{\nu _{2}}+\frac{\lambda }{n_{2}}%
e^{-i\pi \nu _{2}})^{r+1}}\right] dz  \notag \\
&=&\frac{\lambda e^{i\pi \nu _{2}}}{2\pi i}\int_{0}^{\infty }\frac{e^{-zt}dz%
}{n_{1}e^{i\pi (\nu _{2}-\nu _{1})}z^{\nu _{1}}+n_{2}z^{\nu _{2}}+\lambda
e^{i\pi \nu _{2}}}-\frac{\lambda e^{-i\pi \nu _{2}}}{2\pi i}\int_{0}^{\infty
}\frac{e^{-zt}dz}{n_{1}e^{-i\pi (\nu _{2}-\nu _{1})}z^{\nu _{1}}+n_{2}z^{\nu
_{2}}+\lambda e^{-i\pi \nu _{2}}}  \notag \\
&=&\left[ w=\frac{z}{t}\right]  \notag \\
&=&\frac{\lambda }{2\pi i}\int_{0}^{\infty }e^{-w}dw\left[ \frac{e^{i\pi \nu
_{2}}}{n_{1}e^{i\pi (\nu _{2}-\nu _{1})}w^{\nu _{1}}t^{1-\nu
_{1}}+n_{2}w^{\nu _{2}}t^{1-\nu _{2}}+\lambda e^{i\pi \nu _{2}}t}\right. +
\notag \\
&&\left. -\frac{e^{-i\pi \nu _{2}}}{n_{1}e^{-i\pi (\nu _{2}-\nu _{1})}w^{\nu
_{1}}t^{1-\nu _{1}}+n_{2}w^{\nu _{2}}t^{1-\nu _{2}}+\lambda e^{-i\pi \nu
_{2}}t}\right]  \notag \\
&\simeq &\frac{\lambda \sin (\pi \nu _{2})t^{\nu _{2}-1}\Gamma (1-\nu _{2})}{%
\pi n_{2}}=\frac{\lambda t^{\nu _{2}-1}}{n_{2}\Gamma (\nu _{2})},\qquad
t\rightarrow 0  \notag
\end{eqnarray}%
which shows that, for $t\rightarrow 0$, the asymptotic behavior of $%
\widetilde{f}_{1}^{\nu }$ depends only on the larger fractional index $\nu
_{2}$. The same conclusion can be drawn by looking at the series expansion
of $\widetilde{f}_{1}^{\nu }(t)$ given in (\ref{tem2}).

For $t\rightarrow \infty $, from the sixth line of (\ref{inter}), we have
instead that%
\begin{eqnarray}
&&\widetilde{f}_{1}^{\nu }(t)  \label{inter2} \\
&=&\frac{\lambda t^{-1}}{2\pi i}\int_{0}^{\infty }e^{-w}dw\cdot  \notag \\
&&\cdot \frac{n_{1}w^{\nu _{1}}t^{-\nu _{1}}(e^{i\pi \nu _{1}}-e^{-i\pi \nu
_{1}})+n_{2}w^{\nu _{2}}t^{-\nu _{2}}(e^{i\pi \nu _{2}}-e^{-i\pi \nu _{2}})}{%
(n_{1}e^{i\pi (\nu _{2}-\nu _{1})}w^{\nu _{1}}t^{-\nu _{1}}+n_{2}w^{\nu
_{2}}t^{-\nu _{2}}+\lambda e^{i\pi \nu _{2}})\left( n_{1}e^{-i\pi (\nu
_{2}-\nu _{1})}w^{\nu _{1}}t^{-\nu _{1}}+n_{2}w^{\nu _{2}}t^{-\nu
_{2}}+\lambda e^{-i\pi \nu _{2}}\right) }  \notag \\
&\simeq &\frac{n_{1}\nu _{1}t^{-1-\nu _{1}}\sin (\pi \nu _{1})\Gamma (\nu
_{1})}{\lambda \pi }=\frac{n_{1}\nu _{1}t^{-1-\nu _{1}}}{\lambda \Gamma
(1-\nu _{1})},\qquad t\rightarrow \infty ,  \notag
\end{eqnarray}%
which depends only on the smaller fractional index $\nu _{1.}$ Both
asymptotic expressions (\ref{inter}) and (\ref{inter2}) are exactly the same
as for a fractional Poisson process of single order equal to $\nu _{2}$ and $%
\nu _{1}$ respectively (see \cite{BO}, formulae (2.38) and (2.36)).

Analogously we can analyze the asymptotics of the survival probability,
which turns out to be, for $t\rightarrow \infty ,$
\begin{eqnarray}
&&\widetilde{\Psi }^{\nu }(t)  \label{ast} \\
&=&\frac{n_{2}t^{-\nu _{2}}e^{i\pi }}{2\pi i}\int_{0}^{\infty }\frac{%
e^{-w}w^{\nu _{2}-1}dw}{n_{1}e^{i\pi (\nu _{2}-\nu _{1})}w^{\nu _{1}}t^{-\nu
_{1}}+n_{2}w^{\nu _{2}}t^{-\nu _{2}}+\lambda e^{i\pi \nu _{2}}}+  \notag \\
&&-\frac{n_{2}t^{-\nu _{2}}e^{-i\pi }}{2\pi i}\int_{0}^{\infty }\frac{%
e^{-w}w^{\nu _{2}-1}dw}{n_{1}e^{-i\pi (\nu _{2}-\nu _{1})}w^{\nu
_{1}}t^{-\nu _{1}}+n_{2}w^{\nu _{2}}t^{-\nu _{2}}+\lambda e^{-i\pi \nu _{2}}}%
+  \notag \\
&&+\frac{n_{1}t^{-\nu _{1}}e^{i\pi +i\pi (\nu _{2}-\nu _{1})}}{2\pi i}%
\int_{0}^{\infty }\frac{e^{-w}w^{\nu _{1}-1}dw}{n_{1}e^{i\pi (\nu _{2}-\nu
_{1})}w^{\nu _{1}}t^{-\nu _{1}}+n_{2}w^{\nu _{2}}t^{-\nu _{2}}+\lambda
e^{i\pi \nu _{2}}}+  \notag \\
&&-\frac{n_{1}t^{-\nu _{1}}e^{-i\pi -i\pi (\nu _{2}-\nu _{1})}}{2\pi i}%
\int_{0}^{\infty }\frac{e^{-w}w^{\nu _{1}-1}dw}{n_{1}e^{-i\pi (\nu _{2}-\nu
_{1})}w^{\nu _{1}}t^{-\nu _{1}}+n_{2}w^{\nu _{2}}t^{-\nu _{2}}+\lambda
e^{-i\pi \nu _{2}}}  \notag \\
&\simeq &\frac{n_{1}t^{-\nu _{1}}\sin (\pi (1-\nu _{1}))\Gamma (\nu _{1})}{%
\lambda \pi }=\frac{n_{1}t^{-\nu _{1}}}{\lambda \Gamma (1-\nu _{1})},\qquad
t\rightarrow \infty .  \notag
\end{eqnarray}%
For $t\rightarrow 0,$ by writing down the first terms of the series
expansion in (\ref{tem4}) (at least for $j=0,1,2$ and $r=0,1,2$) and doing
some manipulations, we finally get%
\begin{equation}
\widetilde{\Psi }^{\nu }(t)=1-\frac{\lambda }{n_{2}}\frac{t^{\nu _{2}}}{%
\Gamma \left( \nu _{2}+1\right) }+o(t^{\nu _{2}}),\quad 0\leq t<<1.
\label{ast2}
\end{equation}%
As far as the renewal function is concerned, its asymptotic behavior can be
represented as follows:%
\begin{equation}
\widetilde{m}^{\nu }(t)\simeq \frac{\lambda t^{\nu _{1}}}{n_{1}\Gamma (\nu
_{1}+1)},\qquad t\rightarrow \infty  \label{mt}
\end{equation}%
and%
\begin{equation}
\widetilde{m}^{\nu }(t)\simeq \frac{\lambda t^{\nu _{2}}}{n_{2}\Gamma (\nu
_{2}+1)},\qquad t\rightarrow 0.
\end{equation}%
From (\ref{mt}) it is evident that the mean waiting time, which coincides
with $lim_{t\rightarrow \infty }t/\widetilde{m}^{\nu }(t)$, is infinite,
since $\nu _{1}<1.$

\

\noindent \textbf{Remark 2.7 }We remark that, in Theorem 2.6, the survival
probability $\widetilde{\Psi }^{\nu }$ expressed in (\ref{tem4}) is proved
to be a solution of the relaxation equation of distributed order (\ref{tem7}%
) under the double-order hypothesis (\ref{bi3}). This result can be compared
to the analysis in \cite{Vib}, where only the Laplace transform of the
solution is presented, together with its asymptotic behavior. Formulae (\ref%
{ast}) and (\ref{ast2}) coincide with the result (4.13) obtained therein,
but here we provide an explicit formula of the solution, in terms of
infinite sums of GML functions.

\subsection{Interpolation between fractional and integer-order equation}

We analyze now the following equation:%
\begin{equation}
n_{1}\frac{d^{\nu }p_{k}}{dt^{\nu }}+n_{2}\frac{dp_{k}}{dt}=-\lambda
(p_{k}-p_{k-1}),\quad k\geq 0,\text{ }\nu \in \left( 0,1\right)  \label{par1}
\end{equation}%
which is obtained from (\ref{bi2}), as a special case for $\nu _{2}=1$,
under the usual initial conditions%
\begin{equation}
p_{k}(0)=\left\{
\begin{array}{c}
1\qquad k=0 \\
0\qquad k\geq 1%
\end{array}%
\right. ,  \label{par2}
\end{equation}%
and $p_{-1}(t)=0$. Equation (\ref{par1}) represents an interpolation between
the standard and the fractional equation governing the Poisson process.
Hence the solution, which will be denoted in this case by $\widehat{p}%
_{k}^{\nu }$, must coincide, for $n_{1}=0$, $n_{2}=1$ with the distribution
of the homogeneous Poisson process, i.e. $p_{k},k\geq 0$. On the other hand,
for $n_{1}=1$, $n_{2}=0$ we must retrieve the distribution of the fractional
Poisson process, i.e. $p_{k}^{\nu },k\geq 0$, given in (\ref{gml3}).

The Laplace transform of the solution to equation (\ref{par1}) can be
obtained directly by putting $\nu _{1}=\nu $ and $\nu _{2}=1$ in the result
of Theorem 2.1, so that we get%
\begin{equation}
\mathcal{L}\left\{ \widetilde{p}_{k}^{\nu };\eta \right\} =\frac{\lambda
^{k}n_{1}\eta ^{\nu _{1}-1}+\lambda ^{k}n_{2}}{(\lambda +n_{1}\eta ^{\nu
_{1}}+n_{2}\eta )^{k+1}},\qquad k\geq 0.  \label{par4}
\end{equation}%
for any $k\geq 0.$ In order to invert this expression we adapt the result of
Theorem 2.2 as follows.

\

\noindent \textbf{Theorem 2.7} \ \textit{The solution to equation (\ref{par1}%
), under conditions (\ref{par2}), are given, for any }$k\geq 0,$\textit{\
and }$t>0,$\textit{\ by}%
\begin{eqnarray}
&&\widetilde{p}_{k}^{\nu }(t)=\int_{0}^{+\infty }p_{k}(y)q_{\nu }(y,t)dy
\label{par3} \\
&=&\frac{1}{k!}\int_{0}^{+\infty }y^{k}e^{-y}\left[ n_{1}I^{\nu }(\overline{%
\overline{p}}_{\nu }(\cdot ;y))(t)+n_{2}\overline{\overline{p}}_{\nu }(t;y))%
\right] dy.  \notag
\end{eqnarray}%
\textit{Here }$\overline{\overline{p}}_{\nu }(\cdot ;y)$\textit{\ denotes
the stable law of the random variable }$X_{\nu }$\textit{\ of index }$\nu
\in \left( 0,1\right) $\textit{\ and parameters equal to }$\beta =1,$\textit{%
\ }$\mu =n_{2}|y|/\lambda $\textit{\ and }$\sigma =\left( \frac{n_{1}}{%
\lambda }|y|\cos \frac{\pi \nu }{2}\right) ^{1/\nu }.$

\noindent \textbf{Proof \ }We observe that (\ref{par4}) can be rewritten as
follows%
\begin{equation}
\mathcal{L}\left\{ \widetilde{p}_{k}^{\nu };\eta \right\} =(\frac{n_{1}}{%
\lambda }\eta ^{\nu -1}+\frac{n_{2}}{\lambda })\mathcal{L}\left\{ p_{k};%
\frac{n_{1}}{\lambda }\eta ^{\nu }+\frac{n_{2}}{\lambda }\eta \right\} ,
\notag
\end{equation}%
so that we get, by an argument analugous to Theorem 2.2,%
\begin{eqnarray}
\widetilde{p}_{k}^{\nu }(t) &=&\frac{n_{1}}{\lambda \Gamma (1-\nu )}%
\int_{0}^{t}(t-w)^{-\nu }\left( \int_{0}^{+\infty }p_{k}(y)g_{\nu
}(w;y)dy\right) dw+  \label{par6} \\
&&+\frac{n_{2}}{\lambda }\int_{0}^{+\infty }p_{k}(y)g_{\nu }(t;y)dy  \notag
\end{eqnarray}%
The density in (\ref{par6}) can be expressed as follows%
\begin{eqnarray*}
g_{\nu }(w;y) &=&\mathcal{L}^{-1}\left\{ e^{-(\frac{n_{1}}{\lambda }\eta
^{\nu }+\frac{n_{2}}{\lambda }\eta )|y|};w\right\} \\
&=&\int_{0}^{w}\overline{p}_{\nu }(w-x;y)\delta (x-\frac{n_{2}|y|}{\lambda }%
)dx \\
&=&\overline{\overline{p}}_{\nu }(w;y);
\end{eqnarray*}%
hence%
\begin{eqnarray*}
&&\widetilde{p}_{k}^{\nu }(t)=\int_{0}^{+\infty }p_{k}(y)\left[ \frac{n_{1}}{%
\lambda \Gamma (1-\nu )}\int_{0}^{t}(t-w)^{-\nu }\overline{\overline{p}}%
_{\nu }(w;y)dw+\frac{n_{2}}{\lambda }\overline{\overline{p}}_{\nu }(w;y)%
\right] dy \\
&=&\int_{0}^{+\infty }p_{k}(y)\left[ \frac{n_{1}}{\lambda }I^{\nu }(%
\overline{\overline{p}}_{\nu }(w;y))+\frac{n_{2}}{\lambda }\overline{%
\overline{p}}_{\nu }(w;y)\right] dy,
\end{eqnarray*}%
which coincides with (\ref{par3}).\hfill $\square $

\

\noindent \textbf{Remark 2.8 \ }For $n_{2}=0$ and $n_{1}=1$, we get that $%
\overline{\overline{p}}_{\nu }(t;y)=\overline{p}_{\nu }(t;y)$. Therefore
formula (\ref{par3}) reduces to%
\begin{eqnarray}
&&\widetilde{p}_{k}^{\nu }(t)=\frac{1}{k!}\int_{0}^{+\infty
}y^{k}e^{-y}I^{\nu }(\overline{p}_{\nu }(\cdot ;y))(t)dy  \notag \\
&=&\frac{1}{k!}\int_{0}^{+\infty }y^{k}e^{-y}\overline{v}_{2\nu }(\cdot
;y))(t)dy
\end{eqnarray}%
as for the single-order fractional equation (see (\ref{spe}) and (\ref%
{due.16b})-(\ref{due.16f})). On the other hand, for $n_{1}=0$ and $n_{2}=1$,
it is $\overline{\overline{p}}_{\nu }(t;y)=\delta (y)$ and $\widetilde{p}%
_{k}^{\nu }(t)=p_{k}(t),$ since, in this case, equation (\ref{par1}) reduces
to the equation governing the Poisson distribution.

\

\noindent \textbf{Remark 2.9 }A particular feature in this section is that
for the process governed by (\ref{par1}) the probability generating function
$\widetilde{G}_{\nu }$, as well as the probability of zero events $%
\widetilde{p}_{0}^{\nu }$ and the interarrival time density $\widetilde{f}%
_{1}^{\nu }(t),$ can be expressed as infinite sums of the Kummer confluent
hypergeometric function $_{1}F_{1}\left( \alpha ,\beta ;x\right) .$\ The
latter is defined as\textit{\ }%
\begin{equation}
_{1}F_{1}\left( \alpha ;\gamma ;z\right) =\sum_{j=0}^{\infty }\frac{(\alpha
)_{j}}{(\gamma )_{j}}\frac{z^{j}}{j!},\qquad z,\alpha \in \mathbb{C},\text{ }%
\gamma \in \mathbb{C}\backslash \mathbb{Z}_{0}^{-},  \label{con}
\end{equation}%
where $\left( \gamma \right) _{r}=\gamma (\gamma +1)...(\gamma +r-1)$\ (for $%
r=1,2,...,$\ and $\gamma \neq 0$) and $\left( \gamma \right) _{0}=1,$ or, in
integral form, as%
\begin{equation}
_{1}F_{1}\left( \alpha ;\gamma ;z\right) =\frac{\Gamma (\gamma )}{\Gamma
(\alpha )\Gamma (\gamma -\alpha )}\int_{0}^{1}t^{\alpha -1}(1-t)^{\gamma
-\alpha -1}e^{zt}dt,\qquad 0<\mathcal{R}(\alpha )<\mathcal{R}\text{(}\gamma
),  \label{con2}
\end{equation}%
(see \cite{GR}, p.1085). Indeed it is well-known the following relationship
between the GML function with first parameter equal to one and $%
_{1}F_{1}\left( \alpha ,\gamma ;x\right) $:%
\begin{equation*}
E_{1,\gamma }^{\alpha }(z)=\frac{1}{\Gamma (\gamma )}\,_{1}F_{1}\left(
\alpha ;\gamma ;z\right)
\end{equation*}%
(see \cite{kil}, p.62).

By specializing result (\ref{gu}), the probability generating function $%
\widetilde{G}_{\nu }(u,t)$ is equal to%
\begin{eqnarray}
&&\widetilde{G}_{\nu }(u,t)=\sum_{r=0}^{\infty }\frac{\left( -\frac{%
n_{1}t^{1-\nu }}{n_{2}}\right) ^{r}}{\Gamma (r(1-\nu )+1)}\,_{1}F_{1}\left(
r+1;r(1-\nu )+1;-\frac{\lambda (1-u)t}{n_{2}}\right) +  \notag \\
&&-\sum_{r=0}^{\infty }\frac{\left( -\frac{n_{1}t^{1-\nu }}{n_{2}}\right)
^{r+1}}{\Gamma ((r+1)(1-\nu )+1)}\,_{1}F_{1}\left( r+1;(r+1)(1-\nu )+1;-%
\frac{\lambda (1-u)t}{n_{2}}\right) .  \notag
\end{eqnarray}%
Analogously, from (\ref{tem2}) the interarrival time density reads, in this
case,%
\begin{equation}
\widetilde{f}_{1}^{\nu }(t)=\frac{\lambda }{n_{2}}\sum_{r=0}^{\infty }\frac{%
\left( -\frac{n_{1}t^{1-\nu }}{n_{2}}\right) ^{r}}{\Gamma (r-\nu r+1)}%
\,_{1}F_{1}\left( r+1;r-\nu r+1;-\frac{\lambda t}{n_{2}}\right) .
\end{equation}

\

\noindent \textbf{Remark 2.10 }The expected value of the renewal process $%
\widetilde{\mathcal{N}}_{\nu }(t),t>0$ with distribution (\ref{par3}) is
given by%
\begin{equation}
\mathbb{E}\widetilde{\mathcal{N}}_{\nu }(t)=\frac{\lambda t}{n_{2}}E_{1-\nu
,2}\left( -\frac{n_{1}t^{1-\nu _{1}}}{n_{2}}\right) ,  \label{exp3}
\end{equation}%
so that we get this asymptotic behavior%
\begin{equation*}
\mathbb{E}\widetilde{\mathcal{N}}_{\nu }(t)\simeq \frac{\lambda t^{\nu }}{%
n_{1}\Gamma (1+\nu )},\quad t\rightarrow \infty .
\end{equation*}%
This expression coincides with (\ref{exp2}), for $n_{1}=1$: the mean value
is not influenced by the presence of the first derivative. On the contrary,
for $t\rightarrow 0$, we obtain from (\ref{exp3}) that%
\begin{equation*}
\mathbb{E}\widetilde{\mathcal{N}}_{\nu }(t)\simeq \frac{\lambda t}{n_{2}}%
,\quad t\rightarrow 0,
\end{equation*}%
i.e. the usual expected value of the Poisson process. Therefore the first
derivative dominates equation (\ref{par1}) asymptotically, as $t\rightarrow
0.$

\section{Diffusion equations of distributed order\ }

We study equation (\ref{p4}) in the double-order hypothesis (\ref{bi3}), i.e.%
\begin{equation}
n_{1}\frac{\partial ^{\nu _{1}}v}{\partial t^{\nu _{1}}}+n_{2}\frac{\partial
^{\nu _{2}}v}{\partial t^{\nu _{2}}}=\frac{\partial ^{2}v}{\partial x^{2}}%
,\quad x\in \mathbb{R},t>0,\;v(x,0)=\delta (x),\;n_{1},n_{2}>0,  \label{tre2}
\end{equation}%
for $0<\nu _{1}<\nu _{2}\leq 1$. Equation (\ref{tre2}) can be viewed also as
the particular case (for $\gamma =2)$ of equation (\ref{gu3}) analyzed in
\cite{sa2}$.$ In that paper only the Fourier transform of the solution is
given in explicit form, in terms of infinite sums of GML functions. Our aim
is to give an explicit form of the solution, by using an approach similar to
the previous section and providing an expression of the density of the
random time in the subordinating relationship (\ref{sub}). This turns out to
coincide with the density of the random time $\widetilde{\mathcal{T}}_{\nu
_{1},\nu _{2}}(t)$, i.e. with $q_{\nu _{1}.\nu _{2}}$ given in (\ref{bi7})
or (\ref{bi11}).

\

\noindent \textbf{Theorem 3.1} \ \textit{The solution to equation (\ref{tre2}%
), is given by}%
\begin{eqnarray}
&&\widetilde{v}_{\nu _{1},\nu _{2}}(x,t)=\int_{0}^{+\infty }f(x,y)q_{\nu
_{1},\nu _{2}}(y,t)dy  \label{tre3} \\
&=&\int_{0}^{+\infty }\frac{e^{-x^{2}/4y}}{\sqrt{4\pi y}}\left[ \int_{0}^{t}%
\overline{p}_{\nu _{2}}(t-s;y)\overline{v}_{2\nu _{1}}(y,s)ds+\int_{0}^{t}%
\overline{p}_{\nu _{1}}(t-s;y)\overline{v}_{2\nu _{2}}(y,s)ds\right] dy,
\notag
\end{eqnarray}%
\textit{where }$f$\textit{\ is the transition density of a Brownian motion }$%
B(t),t>0,$\textit{\ }$\overline{v}_{2\nu _{j}}$\textit{\ is given in (\ref%
{due.16b})-(\ref{due.16f}) and }$\overline{p}_{\nu _{j}}(\cdot ;y)$\textit{\
denotes the stable law of the random variable }$X_{\nu _{j}}$\textit{\ of
index }$\nu _{j}\in \left( 0,1\right) $\textit{\ with parameters }$\beta =1,$%
\textit{\ }$\mu =0$\textit{\ and }$\sigma =\left( n_{j}|y|\cos \frac{\pi \nu
_{j}}{2}\right) ^{1/\nu _{j}},$\textit{\ for }$j=1,2.$\textit{\
Alternatively the density }$q_{\nu _{1}.\nu _{2}}$\textit{\ in (\ref{tre3})
can be written as in (\ref{bi11}).}

\noindent \textbf{Proof \ }We take the Fourier transform of (\ref{tre2}), so
that we get%
\begin{equation}
\left\{
\begin{array}{l}
n_{1}\frac{\partial ^{\nu _{1}}\widetilde{V}}{\partial t^{\nu _{1}}}+n_{2}%
\frac{\partial ^{\nu _{2}}\widetilde{V}}{\partial t^{\nu _{2}}}=-\theta ^{2}%
\widetilde{V} \\
\widetilde{V}(\theta ,0)=1%
\end{array}%
\right.  \label{tre4}
\end{equation}%
where
\begin{equation*}
\widetilde{V}_{\nu _{1},\nu _{2}}(\theta ,t)=\mathcal{F}\left\{ \widetilde{v}%
_{\nu _{1},\nu _{2}};\theta \right\} =\int_{-\infty }^{+\infty }e^{i\theta x}%
\widetilde{v}_{\nu _{1},\nu _{2}}(x,t)dx.
\end{equation*}%
Taking now the Laplace transform of (\ref{tre4}) we get%
\begin{equation}
\mathcal{L}\left\{ \widetilde{V}_{\nu _{1},\nu _{2}};\theta ,\eta \right\} =%
\frac{n_{1}\eta ^{\nu _{1}-1}+n_{2}\eta ^{\nu _{2}-1}}{n_{1}\eta ^{\nu
_{1}}+n_{2}\eta ^{\nu _{2}}+\theta ^{2}}.  \label{tre6}
\end{equation}%
We can invert the Laplace transform, by noting that it coincides with (\ref%
{lo}) for $\lambda =\theta ^{2}$, as follows:%
\begin{eqnarray}
&&\widetilde{V}_{\nu _{1},\nu _{2}}(\theta ,t)  \label{tre5} \\
&=&\sum_{r=0}^{\infty }\left( -\frac{n_{1}t^{\nu _{2}-\nu _{1}}}{n_{2}}%
\right) ^{r}E_{\nu _{2},(\nu _{2}-\nu _{1})r+1}^{r+1}\left( -\frac{\theta
^{2}t^{\nu _{2}}}{n_{2}}\right) -\sum_{r=0}^{\infty }\left( -\frac{%
n_{1}t^{\nu _{2}-\nu _{1}}}{n_{2}}\right) ^{r+1}E_{\nu _{2},(\nu _{2}-\nu
_{1})(r+1)+1}^{r+1}\left( -\frac{\theta ^{2}t^{\nu _{2}}}{n_{2}}\right) ,
\notag
\end{eqnarray}%
thus obtaining a first form for the solution to (\ref{tre4}). Since
inverting the Fourier transform (\ref{tre5}) seems not possible in closed
form, we rewrite (\ref{tre6}) as follows:%
\begin{equation*}
\mathcal{L}\left\{ \widetilde{V}_{\nu _{1},\nu _{2}};\theta ,\eta \right\}
=(n_{1}\eta ^{\nu _{1}-1}+n_{2}\eta ^{\nu _{2}-1})\frac{1}{2}%
\int_{0}^{+\infty }e^{-n_{1}w\eta ^{\nu _{1}}+n_{2}w\eta ^{\nu _{2}}-\theta
^{2}w}dw.
\end{equation*}%
We note that the term $e^{-(n_{1}\eta ^{\nu _{1}}+n_{2}\eta ^{\nu _{2}})w}$
can be seen again as the convolution of two stable laws $\overline{p}_{\nu
_{j}}$ of index $\nu _{j}\in \left( 0,1\right) $ and parameters equal to $%
\beta =1,$ $\mu =0$ and $\sigma =\left( \frac{1}{2}n_{j}|w|\cos \frac{\pi
\nu _{j}}{2}\right) ^{1/\nu _{j}}$ for $j=1,2$ (see (\ref{sta})). Therefore
we get, alternatively to (\ref{tre5})%
\begin{eqnarray}
&&\widetilde{V}_{\nu _{1},\nu _{2}}(\theta ,t)  \label{tre7} \\
&=&\frac{n_{1}}{\Gamma (1-\nu _{1})}\int_{0}^{t}(t-z)^{-\nu
_{1}}dz\int_{0}^{+\infty }e^{-\theta ^{2}w}\left[ \int_{0}^{z}\overline{p}%
_{\nu _{1}}(z-x;w)\overline{p}_{\nu _{2}}(x;w)dx\right] dw+  \notag \\
&&+\frac{n_{2}}{\Gamma (1-\nu _{2})}\int_{0}^{t}(t-z)^{-\nu
_{2}}dz\int_{0}^{+\infty }e^{-\theta ^{2}w}\left[ \int_{0}^{z}\overline{p}%
_{\nu _{1}}(z-x;w)\overline{p}_{\nu _{2}}(x;w)dx\right] dw  \notag \\
&=&n_{1}\int_{0}^{+\infty }e^{-\theta ^{2}w}\left[ \int_{0}^{t}I^{\nu
_{1}}\left\{ \overline{p}_{\nu _{1}}(\cdot ;w)\right\} (x)\,\overline{p}%
_{\nu _{2}}(t-x;w)dx\right] dw+  \notag \\
&&+n_{2}\int_{0}^{+\infty }e^{-\theta ^{2}w}\left[ \int_{0}^{t}I^{\nu
_{2}}\left\{ \overline{p}_{\nu _{2}}(\cdot ;w)\right\} (x)\,\overline{p}%
_{\nu _{1}}(t-x;w)dx\right] dw  \notag \\
&=&\int_{0}^{+\infty }e^{-\theta ^{2}w}\left[ \int_{0}^{t}\overline{v}_{2\nu
_{1}}(x,w)\,\overline{p}_{\nu _{2}}(t-x;w)dx\right] dw+  \notag \\
&&+\int_{0}^{+\infty }e^{-\theta ^{2}w}\left[ \int_{0}^{t}\overline{v}_{2\nu
_{2}}(x,w)\,\overline{p}_{\nu _{1}}(t-x;w)dx\right] dw,  \notag
\end{eqnarray}%
where again $\overline{v}_{2\nu _{j}}$ is the solution to equation (\ref%
{due.16}) with $c^{2}=1/n_{j},$ $j=1,2.$ Finally, we recognize in (\ref{tre7}%
) the Fourier transform of the Gaussian density, with variance $2|w|,$ so
that we can write the subordination relationship (\ref{tre3}).\hfill $%
\square $

\

The previous theorem shows that the solution to the double-order equation (%
\ref{tre2}) can be seen as the density of the random-time process%
\begin{equation}
\widetilde{\mathcal{B}}_{\nu _{1},\nu _{2}}(t)\equiv B(\widetilde{\mathcal{T}%
}_{\nu _{1},\nu _{2}}(t)),\qquad t>0,  \label{bt}
\end{equation}%
where $B$ is a Brownian motion (with infinitesimal variance equal to $2$)
and $\widetilde{\mathcal{T}}_{\nu _{1},\nu _{2}}$ is the random time,
independent from $B,$ with density $q_{\nu _{1}.\nu _{2}}$ given in (\ref%
{tre3}) or, alternatively, in (\ref{bi11}), for $\lambda =1$. By using the
results obtained in Theorem 2.5, we can evaluate the moments of $\widetilde{%
\mathcal{B}}_{\nu _{1},\nu _{2}}$, as follows:

\

\noindent \textbf{Theorem 3.2 \ }\textit{The }$k$\textit{-th order moments
of the process }$\widetilde{\mathcal{B}}_{\nu _{1},\nu _{2}}$\textit{\ are
given by}

\begin{eqnarray*}
&&\mathbb{E}\widetilde{\mathcal{B}}_{\nu _{1},\nu _{2}}^{k}(t) \\
&=&\left\{
\begin{array}{l}
0\qquad \text{for }k=2h+1 \\
\frac{n_{1}t^{\nu _{2}h+\nu _{2}-\nu _{1}}}{n_{2}^{h+1}}(2h)!E_{\nu _{2}-\nu
_{1},\nu _{2}h+\nu _{2}-\nu _{1}+1}^{h+1}\left( -\frac{n_{1}t^{\nu _{2}-\nu
_{1}}}{n_{2}}\right) +\frac{t^{\nu _{2}h}}{n_{2}^{h}}(2h)!E_{\nu _{2}-\nu
_{1},\nu _{2}h+1}^{h+1}\left( -\frac{n_{1}t^{\nu _{2}-\nu _{1}}}{n_{2}}%
\right) \\
\ \qquad \qquad \qquad \qquad \qquad \qquad \qquad \qquad \qquad \qquad
\qquad \qquad \qquad \text{for }k=2h%
\end{array}%
\right.
\end{eqnarray*}%
\textbf{Proof }By the definition (\ref{bt}) we can write%
\begin{equation*}
\mathbb{E}\widetilde{\mathcal{B}}_{\nu _{1},\nu _{2}}^{k}(t)=\mathbb{E}\left[
B(\widetilde{\mathcal{T}}_{\nu _{1},\nu _{2}}(t))\right] ^{k}=\int_{-\infty
}^{+\infty }x^{k}\int_{0}^{+\infty }f(x,y)q_{\nu _{1},\nu _{2}}(y,t)dydx.
\end{equation*}%
The odd order moments of $\widetilde{\mathcal{B}}_{\nu _{1},\nu _{2}}$ are
obviously null, while the moments of order $2h,$ for $h\in \mathbb{N}$, can
be evaluated as follows:%
\begin{eqnarray*}
&&\mathbb{E}\widetilde{\mathcal{B}}_{\nu _{1},\nu _{2}}^{2h}(t) \\
&=&\int_{0}^{+\infty }q_{\nu _{1},\nu _{2}}(y,t)\int_{-\infty }^{+\infty
}x^{2h}\frac{e^{-x^{2}/4y}}{\sqrt{4\pi y}}dxdy \\
&=&\left[ x=\sqrt{4yw}\right] \\
&=&\frac{1}{\sqrt{\pi }}\int_{0}^{+\infty }q_{\nu _{1},\nu
_{2}}(y,t)\int_{0}^{+\infty }\left( 4yw\right) ^{h}\frac{e^{-w}}{\sqrt{w}}%
dwdy \\
&=&\frac{2^{2h}}{\sqrt{\pi }}\Gamma \left( h+\frac{1}{2}\right)
\int_{0}^{+\infty }y^{h}q_{\nu _{1},\nu _{2}}(y,t)dy \\
&=&2\frac{\left( 2h-1\right) !}{(h-1)!}\int_{0}^{+\infty }y^{h}q_{\nu
_{1},\nu _{2}}(y,t)dy \\
&=&\frac{n_{1}t^{\nu _{2}h+\nu _{2}-\nu _{1}}}{n_{2}^{h+1}}(2h)!E_{\nu
_{2}-\nu _{1},\nu _{2}h+\nu _{2}-\nu _{1}+1}^{h+1}\left( -\frac{n_{1}t^{\nu
_{2}-\nu _{1}}}{n_{2}}\right) +\frac{t^{\nu _{2}h}}{n_{2}^{h}}(2h)!E_{\nu
_{2}-\nu _{1},\nu _{2}h+1}^{h+1}\left( -\frac{n_{1}t^{\nu _{2}-\nu _{1}}}{%
n_{2}}\right) ,
\end{eqnarray*}%
where, in the last step, we have applied formula (\ref{fac}) and the
relationship%
\begin{equation}
\mathbb{E}\left[ \widetilde{\mathcal{T}}_{\nu _{1},\nu _{2}}(t)\right] ^{k}=%
\mathbb{E}\left[ \widetilde{\mathcal{N}}_{\nu _{1},\nu _{2}}(t)\left(
\widetilde{\mathcal{N}}_{\nu _{1},\nu _{2}}(t)-1\right) ...(\widetilde{%
\mathcal{N}}_{\nu _{1},\nu _{2}}(t)-k+1)\right] ,\text{\qquad }k\in \mathbb{N%
}  \label{rel}
\end{equation}%
proved in Theorem 2.5.\hfill $\square $

\

\noindent \textbf{Remark 3.1 \ }We can check the previous result by noting
that, for $h=1$, we get the second-order moment obtained in \cite{che1} (see
formula (16), for $\tau =D=1$):%
\begin{eqnarray}
&&\mathbb{E}\widetilde{\mathcal{B}}_{\nu _{1},\nu _{2}}^{2}(t)  \label{e7} \\
&=&\frac{2n_{1}t^{2\nu _{2}-\nu _{1}}}{n_{2}^{2}}\sum_{j=0}^{\infty }\frac{%
(j+1)!\left( -\frac{n_{1}t^{\nu _{2}-\nu _{1}}}{n_{2}}\right) ^{j}}{j!\Gamma
\left( \left( \nu _{2}-\nu _{1}\right) j+2\nu _{2}-\nu _{1}+1\right) }+\frac{%
2t^{\nu _{2}-\nu _{1}}}{n_{2}}\sum_{j=0}^{\infty }\frac{(j+1)!\left( -\frac{%
n_{1}t^{\nu _{2}-\nu _{1}}}{n_{2}}\right) ^{j}}{j!\Gamma \left( \left( \nu
_{2}-\nu _{1}\right) j+\nu _{2}+1\right) }  \notag \\
&=&-\frac{2t^{\nu _{2}}}{n_{2}}\sum_{l=0}^{\infty }\frac{l\left( -\frac{%
n_{1}t^{\nu _{2}-\nu _{1}}}{n_{2}}\right) ^{l}}{\Gamma \left( \left( \nu
_{2}-\nu _{1}\right) l+\nu _{2}+1\right) }+\frac{2t^{\nu _{2}}}{n_{2}}%
\sum_{j=0}^{\infty }\frac{(j+1)\left( -\frac{n_{1}t^{\nu _{2}-\nu _{1}}}{%
n_{2}}\right) ^{j}}{\Gamma \left( \left( \nu _{2}-\nu _{1}\right) j+\nu
_{2}+1\right) }  \notag \\
&=&\frac{2t^{\nu _{2}}}{n_{2}}E_{\nu _{2}-\nu _{1},\nu _{2}+1}\left( -\frac{%
n_{1}t^{\nu _{2}-\nu _{1}}}{n_{2}}\right) .  \notag
\end{eqnarray}

\

Our attention is now addressed to the equation solved by the density $q_{\nu
_{1},\nu _{2}}$ of the time argument $\widetilde{\mathcal{T}}_{\nu _{1},\nu
_{2}}$ (which is shared by the processes $\widetilde{\mathcal{N}}_{\nu
_{1},\nu _{2}}$ and $\widetilde{\mathcal{B}}_{\nu _{1},\nu _{2}}$ ). In
analogy with the single-order fractional case this equation must be of
\textquotedblleft second-order\textquotedblright\ (involving the two
fractional indexes $\nu _{1},\nu _{2})$, but is not evidently given by (\ref%
{tre2}). We prove in the next theorem that a further time-fractional
derivative must be included in the diffusion equation (\ref{tre2}) in order
to obtain $q_{\nu _{1},\nu _{2}}$ as solution.

\

\noindent \textbf{Theorem 3.3 \ }\textit{The density }$q_{\nu _{1},\nu
_{2}}(x,t)$ \textit{coincides with the folded solution}%
\begin{equation}
\overline{v}(x,t)=\left\{
\begin{array}{l}
2v(x,t),\qquad x\geq 0 \\
0,\qquad x<0%
\end{array}%
\right.  \label{ops}
\end{equation}%
\textit{of the following equation}

\begin{equation}
\left( n_{1}\frac{\partial ^{\nu _{1}}v}{\partial t^{\nu _{1}}}+n_{2}\frac{%
\partial ^{\nu _{2}}v}{\partial t^{\nu _{2}}}\right) ^{2}=\frac{\partial
^{2}v}{\partial x^{2}},\quad x\in \mathbb{R},t>0,\;n_{1},n_{2}>0,  \label{e3}
\end{equation}%
\textit{for} $0<\nu _{1}<\nu _{2}\leq 1$,\textit{\ with initial conditions}%
\begin{equation}
\left\{
\begin{array}{l}
v(x,0)=\delta (x),\text{ for }0<\nu _{1}<\nu _{2}\leq 1 \\
\left. \frac{\partial }{\partial t}v(x,t)\right\vert _{t=0}=0\text{ for }%
\frac{1}{2}<\nu _{1}<\nu _{2}\leq 1%
\end{array}%
\right. .  \label{e5}
\end{equation}

\noindent \textbf{Proof }Let us take the Fourier transform of (\ref{e2}),
which reads

\begin{eqnarray}
\mathcal{L}\left\{ Q_{\nu _{1}.\nu _{2}};\theta ,\eta \right\} &\equiv
&\int_{-\infty }^{\infty }e^{i\theta y}\mathcal{L}\left\{ q_{\nu _{1}.\nu
_{2}};\eta \right\} dy  \label{e6} \\
&=&\left( n_{1}\eta ^{\nu _{1}-1}+n_{2}\eta ^{\nu _{2}-1}\right)
\int_{-\infty }^{\infty }e^{i\theta y}e^{-\left( n_{1}\eta ^{\nu
_{1}}+n_{2}\eta ^{\nu _{2}}\right) |y|}dy  \notag \\
&=&\left( n_{1}\eta ^{\nu _{1}-1}+n_{2}\eta ^{\nu _{2}-1}\right) \left[
\frac{1}{n_{1}\eta ^{\nu _{1}}+n_{2}\eta ^{\nu _{2}}-i\theta }+\frac{1}{%
n_{1}\eta ^{\nu _{1}}+n_{2}\eta ^{\nu _{2}}+i\theta }\right]  \notag \\
&=&\frac{2\left( n_{1}\eta ^{\nu _{1}}+n_{2}\eta ^{\nu _{2}}\right) ^{2}}{%
\eta \left( n_{1}^{2}\eta ^{2\nu _{1}}+n_{2}^{2}\eta ^{2\nu
_{2}}+2n_{1}n_{2}\eta ^{\nu _{1}+\nu _{2}}+\theta ^{2}\right) },  \notag
\end{eqnarray}%
by taking, for simplicity, $\lambda =1.$ We take now the Laplace-Fourier
transform of equation (\ref{e3}), by considering formula (\ref{jen})
together with the initial conditions (\ref{e5}):%
\begin{gather*}
n_{1}^{2}\eta ^{2\nu _{1}}\mathcal{L}\left\{ Q_{\nu _{1}.\nu _{2}};\theta
,\eta \right\} -n_{1}^{2}\eta ^{2\nu _{1}-1}+n_{2}^{2}\eta ^{2\nu _{2}}%
\mathcal{L}\left\{ Q_{\nu _{1}.\nu _{2}};\theta ,\eta \right\}
-n_{2}^{2}\eta ^{2\nu _{2}-1}+ \\
+2n_{1}n_{2}\eta ^{\nu _{1}+\nu _{2}}\mathcal{L}\left\{ Q_{\nu _{1}.\nu
_{2}};\theta ,\eta \right\} -2n_{1}n_{2}\eta ^{\nu _{1}+\nu _{2}-1}=-\theta
^{2}\mathcal{L}\left\{ Q_{\nu _{1}.\nu _{2}};\theta ,\eta \right\} ,
\end{gather*}%
whose solution coincides with (\ref{e6}), by taking into account (\ref{ops}%
).\hfill $\square $

\

\noindent \textbf{Remark 3.2 \ }Equation (\ref{e3}) can be itself
interpreted as a distributed order fractional equation, by assuming a
different density for the random fractional index $\nu ,$ which, in this
case, is defined on the interval $\left( 0,2\right] $: indeed we can
formulate $n(\nu )$, as follows:%
\begin{equation}
n(\nu )=n_{1}^{2}\delta (\nu -2\nu _{1})+n_{2}^{2}\delta (\nu -2\nu
_{2})+2n_{1}n_{2}\delta (\nu -(\nu _{1}+\nu _{2})),\qquad 0<\nu _{1}<\nu
_{2}\leq 1,
\end{equation}%
for $n_{1},n_{2}\geq 0$ and such that $n_{1}+n_{2}=1.$ The last condition is
enough to fulfill the normalization requirement $%
n_{1}^{2}+n_{2}^{2}+2n_{1}n_{2}=1.$

By considering (\ref{fac}) together with (\ref{rel}) we can obtain the
second-order moment of the diffusion process $\widetilde{\mathcal{T}}_{\nu
_{1},\nu _{2}}$ governed by equation (\ref{e3}), i.e.

\begin{eqnarray}
&&\mathbb{E}\left[ \widetilde{\mathcal{T}}_{\nu _{1},\nu _{2}}(t)\right]
^{2}=\frac{2n_{1}t^{3\nu _{2}-\nu _{1}}}{n_{2}^{3}}E_{\nu _{2}-\nu _{1},3\nu
_{2}-\nu _{1}+1}^{3}\left( -\frac{n_{1}t^{\nu _{2}-\nu _{1}}}{n_{2}}\right) +%
\frac{2t^{2\nu _{2}}}{n_{2}^{2}}E_{\nu _{2}-\nu _{1},2\nu _{2}+1}^{3}\left( -%
\frac{n_{1}t^{\nu _{2}-\nu _{1}}}{n_{2}}\right)  \label{e8} \\
&=&\frac{n_{1}t^{3\nu _{2}-\nu _{1}}}{n_{2}^{3}}\sum_{j=0}^{\infty }\frac{%
(j+2)(j+1)}{\Gamma ((\nu _{2}-\nu _{1})j+3\nu _{2}-\nu _{1}+1)}\left( -\frac{%
n_{1}t^{\nu _{2}-\nu _{1}}}{n_{2}}\right) ^{j}+  \notag \\
&&+\frac{t^{2\nu _{2}}}{n_{2}^{2}}\sum_{j=0}^{\infty }\frac{(j+2)(j+1)}{%
\Gamma ((\nu _{2}-\nu _{1})j+2\nu _{2}+1)}\left( -\frac{n_{1}t^{\nu _{2}-\nu
_{1}}}{n_{2}}\right) ^{j}  \notag \\
&=&\left[ l=j+1\right]  \notag \\
&=&-\frac{t^{2\nu _{2}}}{n_{2}^{2}}\sum_{l=1}^{\infty }\frac{(l+1)l}{\Gamma
((\nu _{2}-\nu _{1})l+2\nu _{2}+1)}\left( -\frac{n_{1}t^{\nu _{2}-\nu _{1}}}{%
n_{2}}\right) ^{l}+\frac{t^{2\nu _{2}}}{n_{2}^{2}}\sum_{l=0}^{\infty }\frac{%
(l+2)(l+1)}{\Gamma ((\nu _{2}-\nu _{1})l+2\nu _{2}+1)}\left( -\frac{%
n_{1}t^{\nu _{2}-\nu _{1}}}{n_{2}}\right) ^{l}  \notag \\
&=&\frac{2t^{2\nu _{2}}}{n_{2}^{2}}E_{\nu _{2}-\nu _{1},2\nu
_{2}+1}^{2}\left( -\frac{n_{1}t^{\nu _{2}-\nu _{1}}}{n_{2}}\right) .  \notag
\end{eqnarray}%
We compare (\ref{e8}) with the second-order moment of the diffusion process $%
\widetilde{\mathcal{B}}_{\nu _{1},\nu _{2}}$ governed by equation (\ref{tre2}%
) (which is given in (\ref{e7})), noting that, apart from the similar
structure, the second-order moment of $\widetilde{\mathcal{T}}_{\nu _{1},\nu
_{2}}$ is expressed in terms of a GML function, instead of a standard one.
For $t\rightarrow \infty $ the different asymptotic behavior is described
below, by applying (\ref{asy3}) and (\ref{asy4}):%
\begin{eqnarray}
\mathbb{E}\left[ \widetilde{\mathcal{T}}_{\nu _{1},\nu _{2}}(t)\right] ^{2}
&\simeq &\frac{2t^{2\nu _{1}}}{n_{1}^{2}\Gamma (1+2\nu _{1})}  \label{asy5}
\\
\mathbb{E}\widetilde{\mathcal{B}}_{\nu _{1},\nu _{2}}^{2}(t) &\simeq &\frac{%
2t^{\nu _{1}}}{n_{1}\Gamma (1+\nu _{1})}  \notag
\end{eqnarray}

For $t\rightarrow 0$, since the behavior of the GML and the standard
Mittag-Leffler function coincides, we must apply, in both cases, formula (%
\ref{asy2}). Hence we get:%
\begin{eqnarray}
\mathbb{E}\left[ \widetilde{\mathcal{T}}_{\nu _{1},\nu _{2}}(t)\right] ^{2}
&\simeq &\frac{2t^{2\nu _{2}}}{n_{2}^{2}\Gamma (1+2\nu _{2})}  \label{asy6}
\\
\mathbb{E}\widetilde{\mathcal{B}}_{\nu _{1},\nu _{2}}^{2}(t) &\simeq &\frac{%
2t^{\nu _{2}}}{n_{2}\Gamma (1+\nu _{2})}.  \notag
\end{eqnarray}%
We can conclude from (\ref{asy5}) and (\ref{asy6}) that, in the case where $%
\nu _{1},\nu _{2}<1/2,$ the effect of diffusion with retardation, which is
characteristic of $\widetilde{\mathcal{B}}_{\nu _{1},\nu _{2}}$ (see also
\cite{che1}), is emphasized for the process $\widetilde{\mathcal{T}}_{\nu
_{1},\nu _{2}}$, since the difference between $2\nu _{2}$ and $2\nu _{1}$ is
greater than for $\nu _{2}$ and $\nu _{1}.$

A different conclusion should be drawn in the case where either $\nu _{2}$
or both $\nu _{1},\nu _{2}$ are greater than $1/2.$ For $\nu _{1}<1/2$ and $%
\nu _{2}>1/2$, we can observe that the asymptotic behavior of $\mathbb{E}%
\left[ \widetilde{\mathcal{T}}_{\nu _{1},\nu _{2}}(t)\right] ^{2}$
drastically changes, at least for $t\rightarrow 0$: indeed in this case it
goes to zero faster than $t.$ For $\nu _{1},\nu _{2}>1/2,$ in addition to
this effect, we note that, also for $t\rightarrow \infty ,$ the rate
convergence is greater than in the standard diffusion case (besides that of
diffusion with retardation). Therefore the process governed by (\ref{e3}),
for $\nu _{1},\nu _{2}>1/2$, can be interpreted as a diffusion with
acceleration.

\

\end{document}